\documentclass{amsart}
\usepackage[english]{babel}
\usepackage[latin1]{inputenc}

\usepackage[mathcal]{eucal}
\usepackage{amsmath}
\usepackage{amsfonts}
\usepackage{amssymb}
\usepackage{amsthm}
\usepackage{graphicx,epsfig}
\usepackage{amscd}
\DeclareMathOperator{\re}{Re}

\DeclareMathOperator{\diag}{diag}

\newcommand{\M}{\mathbb{M}}
\newcommand{\R}{\mathbb{R}}
\newcommand{\h}{\mathbb{H}}

\newcommand{\s}{\mathbb{S}}
\newcommand{\LL}{\mathbb{L}}
\newcommand{\E}{\mathbb{E}}

\newcommand{\rmQ}{\mathrm{Q}}
\newcommand{\rmJ}{\mathrm{J}}

\newcommand{\rmT}{\mathrm{T}}
\newcommand{\rmS}{\mathrm{S}}
\newcommand{\rmd}{\mathrm{d}}

\newcommand{\rmR}{\mathrm{R}}

\newcommand{\cU}{{\mathcal U}}
\newcommand{\cM}{{\mathcal M}}

\newcommand{\cF}{{\mathcal F}}

\newcommand{\cV}{{\mathcal V}}

\newcommand{\cH}{{\mathcal H}}
\newcommand{\cA}{{\mathcal A}}

\newcommand{\cC}{{\mathcal C}}
\newcommand{\cD}{{\mathcal D}}

\newcommand{\cZ}{{\mathcal Z}}

\newcommand{\cG}{{\mathcal G}}

\newcommand{\tra}{{}^\mathrm{t}\!}

\begin{document}

\newtheorem{thm}{Theorem}[section]
\newtheorem*{thmintro}{Theorem}
\newtheorem{cor}[thm]{Corollary}
\newtheorem{prop}[thm]{Proposition}
\newtheorem{app}[thm]{Application}
\newtheorem{lemma}[thm]{Lemma}
\newtheorem{ex}[thm]{Example}
\newtheorem{notation}[thm]{Notations}
\newtheorem{hypothesis}[thm]{Hypothesis}

\newtheorem{defin}[thm]{Definition}
\newenvironment{defn}{\begin{defin} \rm}{\end{defin}}
\newtheorem{remk}[thm]{Remark}
\newenvironment{rem}{\begin{remk} \rm}{\end{remk}}

\title[Isometric immersions into $\s^n\times\R$
and $\h^n\times\R$]{Isometric immersions into $\s^n\times\R$
and $\h^n\times\R$ \\ 
and applications to minimal surfaces}
\author{Beno\^\i t Daniel}
\date{}


\address{Universit\'e Paris 7, Institut de Math\'ematiques de Jussieu,
\'Equipe G\'eom\'etrie et Dynamique, Case 7012, 2 place Jussieu, 
75251 Paris Cedex 05, FRANCE}
\email{daniel@math.jussieu.fr}

\begin{abstract}
We give a necessary and sufficient condition for an $n$-dimensional
Riemannian manifold to be isometrically immersed in $\s^n\times\R$ or
$\h^n\times\R$ in terms of its first and second fundamental forms and
of the projection of the vertical vector field on its tangent
plane. We deduce the existence of a one-parameter family of isometric
minimal deformations of a given minimal surface in $\s^2\times\R$ or
$\h^2\times\R$, obtained by rotating the shape operator.
\end{abstract}

\maketitle

\section{Introduction}

It is well known that the first and second fundamental forms of a
hypersurface of a Riemannian manifold satisfy two compatibility
equations called the Gauss and Codazzi equations. More precisely, 
let $\bar\cV$ be an orientable
Riemannian manifold of dimension $n+1$ and $\cV$
a submanifold of $\bar\cV$ of dimension $n$.
Let $\nabla$ (respectively, $\bar\nabla$) be the Riemannian
connection of $\cV$ (respectively, $\bar\cV$), 
$\rmR$ (respectively, $\bar\rmR$) the Riemann curvature tensor
of $\cV$ (respectively, $\bar\cV$), i.e.,
$$\rmR(X,Y)Z=\nabla_Y\nabla_XZ-\nabla_X\nabla_YZ+\nabla_{[X,Y]}Z,$$
and $\rmS$ the shape operator of $\cV$ associated to its unit
normal $N$, i.e.,
$\rmS X=-\bar\nabla_XN$. 
Then the following equations hold for all vector fields
$X,Y,Z$ on $\cV$: 
\begin{equation*}
\rmR(X,Y)Z-\bar\rmR(X,Y)Z
=\langle\rmS X,Z\rangle\rmS Y
-\langle\rmS Y,Z\rangle\rmS X,
\end{equation*}
\begin{equation*}
\nabla_X\rmS Y-\nabla_Y\rmS X-\rmS[X,Y]=\bar\rmR(X,Y)N.
\end{equation*}
These are respectively the Gauss and Codazzi equations.

In the case where $\bar\cV$ is a space form, i.e., the sphere
$\s^{n+1}$, the Euclidean space $\R^{n+1}$ or the hyperbolic space
$\h^{n+1}$, these equations become the following:
\begin{equation} \label{gaussspaceform}
\begin{array}{c}
\langle\rmR(X,Y)Z,W\rangle-\kappa(\langle X,Z\rangle\langle Y,W\rangle
-\langle Y,Z\rangle\langle X,W\rangle) \\
=\langle\rmS X,Z\rangle\langle\rmS Y,W\rangle
-\langle\rmS X,W\rangle\langle\rmS Y,Z\rangle,
\end{array}
\end{equation}
\begin{equation} \label{codazzispaceform}
\nabla_X\rmS Y-\nabla_Y\rmS X-\rmS[X,Y]=0,
\end{equation}
where $\kappa$ is the sectional curvature of $\bar\cV$, i.e.,
$\kappa=1,0,-1$ for $\s^{n+1}$, $\R^{n+1}$ and $\h^{n+1}$
respectively. Thus the Gauss and Codazzi equations only involve the
first and second fundamental forms of $\cV$; they are
defined {\it intrinsicly} on $\cV$ (as soon as we know $\rmS$). This comes
from the fact that these ambiant spaces are isotropic. Moreover, in
this case the Gauss and Codazzi equations are also sufficient
conditions for an 
$n$-dimensional simply connected manifold to be immersed into
$\bar\cV$ with given first and second fundamental forms:
if $\cV$ is a Riemannian manifold endowed with a field $\rmS$ of
symmetric operators $\rmS_y:\rmT_y\cV\to\rmT_y\cV$
such that \eqref{gaussspaceform} and \eqref{codazzispaceform} hold
(where $\rmR$ denotes the Riemann curvature tensor of $\cV$),
then there exists an isometric immersion from $\cV$ into $\bar\cV$
with $\rmS$ as shape operator. The reader can refer to \cite{docarmo},
and also to \cite{tenenblat} for a proof in the case of $\R^{n+1}$.

In the case of a general manifold $\bar\cV$, the Gauss and Codazzi
equations are not defined intrinsicly on $\cV$, since the Riemann
curvature tensor of the ambiant space $\bar\cV$ is involved. Yet, in
the case where $\bar\cV=\s^n\times\R$ or 
$\bar\cV=\h^n\times\R$, these equations are well defined as soon as we
know:
\begin{enumerate}
\item the projection $T$ of the vertical vector 
$\frac{\partial}{\partial t}$ (corresponding to the factor $\R$) onto
the tangent space of $\cV$,
\item the normal component $\nu$ of $\frac{\partial}{\partial t}$, i.e., 
$\nu=\langle N,\frac{\partial}{\partial t}\rangle$.
\end{enumerate}
Indeed, the Gauss and Codazzi equations become the following:
\begin{equation*}
\begin{array}{lll}
\rmR(X,Y)Z & = &
\langle\rmS X,Z\rangle\rmS Y
-\langle\rmS Y,Z\rangle\rmS X \\
& & +\kappa(\langle X,Z\rangle Y-
\langle Y,Z\rangle X
-\langle Y,T\rangle\langle X,Z\rangle T  \\
& & -\langle X,T\rangle\langle Z,T\rangle Y
+\langle X,T\rangle\langle Y,Z\rangle T
+\langle Y,T\rangle\langle Z,T\rangle X),
\end{array}
\end{equation*}
\begin{equation*}
\nabla_X\rmS Y-\nabla_Y\rmS X-\rmS[X,Y]=
\kappa\nu(\langle Y,T\rangle X-\langle X,T\rangle Y),
\end{equation*}
where $\kappa=1$ and $\kappa=-1$ for $\s^n\times\R$ and $\h^n\times\R$
respectively.


The Gauss equation can be formulated in the following equivalent way:
the sectional curvature $K(P)$ (for the metric of $\cV$) of every plane
$P\subset\rmT\cV$ satisfies
$$K(P)=\det\rmS_P+\kappa(1-||T_P||^2)$$ where $\rmS_P$ is
the restriction of $\rmS$ on $P$ and $T_P$ the orthogonal
projection of $T$ on $P$. 

The first aim of this paper is to give a necessary and sufficient
condition in order that a Riemannian manifold with a symmetric
operator $\rmS$ can be isometrically immersed into $\s^n\times\R$ or
$\h^n\times\R$ with $\rmS$ as shape operator. More precisely, we prove
the following theorem. 

\begin{thmintro}[theorem \ref{isometry}]
Let $\cV$ be a simply connected Riemannian manifold of dimension $n$,
$\rmd s^2$ its metric (which we also denote by
$\langle\cdot,\cdot\rangle$) and $\nabla$ its  
Riemannian connection. Let $\rmS$ be a field of symmetric operators
$\rmS_y:\rmT_y\cV\to\rmT_y\cV$, $T$ a vector field on $\cV$
and $\nu$ a smooth function on $\cV$ such that
$||T||^2+\nu^2=1$. 

Let $\M^n=\s^n$ or $\M^n=\h^n$.
Assume that $(\rmd s^2,\rmS,T,\nu)$ satisfies the
Gauss and Codazzi equations for $\M^n\times\R$ and the following
equations:
$$\nabla_XT=\nu\rmS X,\quad
\rmd\nu(X)=-\langle\rmS X,T\rangle.$$
Then there exists an isometric immersion $f:\cV\to\M^n\times\R$ such
that 
the shape operator with respect to the normal $N$ associated to $f$ is
$$\rmd f\circ\rmS\circ\rmd f^{-1}$$ and such that
$$\frac{\partial}{\partial t}=\rmd f(T)+\nu N.$$ Moreover the
immersion is unique up to a global isometry of $\M^n\times\R$
preserving the orientations of both $\M^n$ and $\R$.
\end{thmintro}

The two additional conditions come from the fact that the vertical
vector field $\frac\partial{\partial t}$ is parallel.

The method to prove this theorem is similar to that of Tenenblat
(\cite{tenenblat}): it is based on differential forms, moving frames
and integrable distributions.

This work was motivated by the study of minimal surfaces in
$\s^2\times\R$ and $\h^2\times\R$. There were many recent developments
in the theory of these surfaces. Rosenberg (\cite{rosenbergillinois})
studied the geometry of minimal surfaces in $\s^2\times\R$, and more
generally in $M\times\R$ where $M$ is a surface of non-negative
curvature. Nelli and Rosenberg (\cite{nelli}) studied minimal
surfaces in 
$\h^2\times\R$ and proved a Jenkins-Serrin theorem. Hauswirth
(\cite{hauswirth}) constructed many examples in $\h^2\times\R$. Meeks
and Rosenberg (\cite{meeksrosenberg}) initiated the theory of minimal
surfaces in $M\times\R$ where $M$ is a compact surface.
Recently, Abresch and Rosenberg (\cite{abresch}) extended the notion
of holomorphic Hopf differential to constant mean curvature surfaces
in $\s^2\times\R$ and $\h^2\times\R$; using this holomorphic
differential, they proved that all immersed constant mean curvature
spheres are embedded and rotational.

In this paper, we use our theorem \ref{isometry} to prove the
existence of a one-parameter family of 
isometric minimal deformations of a given minimal surface in
$\s^2\times\R$ or $\h^2\times\R$. This family is obtained by rotating
the shape operator; hence it is the analog of the
associate family of a minimal surface in $\R^3$. This is the
following theorem.

\begin{thmintro}[theorem \ref{associatefamily}]
Let $\Sigma$ be a simply connected Riemann surface and
$x:\Sigma\to\M^2\times\R$ a conformal minimal immersion. Let $N$
be the induced normal. Let
$\rmS$ be the symmetric operator on $\Sigma$ induced by the shape
operator of $x(\Sigma)$. Let $T$
be the vector field on $\Sigma$ such that $\rmd x(T)$ is the
projection of $\frac{\partial}{\partial t}$ onto $\rmT(x(\Sigma))$. Let
$\nu=\left\langle N,\frac{\partial}{\partial t}\right\rangle$.

Let $z_0\in\Sigma$.
Then there exists a unique family $(x_\theta)_{\theta\in\R}$ of
conformal minimal immersions 
$x_\theta:\Sigma\to\M^2\times\R$ such that:
\begin{enumerate}
\item $x_\theta(z_0)=x(z_0)$ and $(\rmd x_\theta)_{z_0}=(\rmd
x)_{z_0}$, 
\item the metrics induced on $\Sigma$ by $x$ and $x_\theta$ are the same,
\item the symmetric operator on $\Sigma$ induced by the shape operator
of $x_\theta(\Sigma)$ is $e^{\theta\rmJ}\rmS$,
\item $\frac{\partial}{\partial t}=\rmd x_\theta(e^{\theta\rmJ}T)+\nu
N_\theta$ where $N_\theta$ is the unit normal to $x_\theta$.
\end{enumerate}

Moreover we have $x_0=x$ and the family $(x_\theta)$ is continuous
with respect to $\theta$.
\end{thmintro}

In particular taking $\theta=\frac{\pi}{2}$ defines a conjugate
surface; the geometric properties of conjugate surfaces in
$\M^2\times\R$ and in $\R^3$ are similar. Finally, we give examples of
conjugate surfaces. In $\s^2\times\R$, we 
show that helicoids and unduloids are conjugate. In $\h^2\times\R$,
we show that helicoids are conjugated to catenoids or to minimal
surfaces foliated by horizontal curves of constant curvature belonging
to the Hauswirth family (see \cite{hauswirth}). 

\section{Preliminaries}

\subsubsection*{Notations.}
In this paper we will use the following index conventions: Latin
letters $i$, $j$, etc, denote integers between $1$ and $n$, Greek
letters $\alpha$, $\beta$, etc, denote integers between $0$ and
$n+1$. For example, the notation $A^i_j=B^i_j$ means that this
relation holds for all integers $i$, $j$ between $1$ and $n$, the
notation $\sum_\alpha C_\alpha$ means $C_0+C_1+\dots+C_{n+1}$.  

The set of vector fields on a Riemannian manifold $\cV$ will be
denoted by $\mathfrak{X}(\cV)$.

We denote by $\frac{\partial}{\partial t}$ the unit vector giving the
orientation of $\R$ in $\M^n\times\R$; we call it the vertical
vector. 

\subsection{The compatibility equations in $\M^n\times\R$}
\label{compatibilitymtimesr}

Let $\M^n=\s^n$ or $\M^n=\h^n$; in the first case we set $\kappa=1$
and in the second case we set $\kappa=-1$. Let $\bar\rmR$ be the
Riemann curvature tensor of $\M^n\times\R$.
Let $\cV$ be an oriented hypersurface of $\M^n\times\R$ and $N$
the unit normal to $\cV$.

\begin{prop}
For $X,Y,Z,W\in\mathfrak{X}(\cV)$ we have
\begin{eqnarray*}
\langle\bar\rmR(X,Y)Z,W\rangle & = &
\kappa(\langle X,Z\rangle\langle Y,W\rangle-
\langle Y,Z\rangle\langle X,W\rangle \\
& &
-\langle Y,T\rangle\langle W,T\rangle\langle X,Z\rangle 
-\langle X,T\rangle\langle Z,T\rangle\langle Y,W\rangle \\
& &
+\langle X,T\rangle\langle W,T\rangle\langle Y,Z\rangle
+\langle Y,T\rangle\langle Z,T\rangle\langle X,W\rangle),
\end{eqnarray*}
\begin{equation*}
\langle\bar\rmR(X,Y)N,Z\rangle=
\kappa\nu
(\langle X,Z\rangle\langle Y,T\rangle-
\langle Y,Z\rangle\langle X,T\rangle),
\end{equation*}
where $$\nu=\left\langle N,\frac{\partial}{\partial t}\right\rangle$$
and $T$ is the projection of $\frac{\partial}{\partial t}$ on $\rmT
\cV$, i.e., $$T=\frac{\partial}{\partial t}-\nu N.$$
\end{prop}

\begin{proof}
Any vector field on $\M^n\times\R$ can be written
$X(m,t)=(X_{\M^n}^t(m),X_{\R}^m(t))$ where, for each $t\in\R$,
$X_{\M^n}^t$ is a vector field on $\M^n$, and, for each $m\in\M^n$,
$X_{\R}^m$ is a vector field on $\R$.
Then for $X,Y,Z,W\in\mathfrak{X}(\M^n\times\R)$ we have
\begin{eqnarray*}
\langle\bar\rmR(X,Y)Z,W\rangle & = &
\langle\bar\rmR_{\M^n}(X_{\M^n},Y_{\M^n})Z_{\M^n},W_{\M^n}\rangle \\
& = &
\kappa(\langle X_{\M^n},Z_{\M^n}\rangle
\langle Y_{\M^n},W_{\M^n}\rangle
-\langle Y_{\M^n},Z_{\M^n}\rangle
\langle X_{\M^n},W_{\M^n}\rangle).
\end{eqnarray*}

We have
$X_{\M^n}=X-\left\langle X,\frac{\partial}{\partial t}\right\rangle
\frac{\partial}{\partial t}$. Thus, if $X\in\rmT\cV$, we have
$X_{\M^n}=X-\langle X,T\rangle
\frac{\partial}{\partial t}$, and similar expressions for
$Y,Z,W\in\rmT\cV$. A computation gives the expected formula for
$\langle\bar\rmR(X,Y)Z,W\rangle$.

Finally we have $N_{\M^n}=N-\nu\frac{\partial}{\partial t}$, so a
computation gives the expected formula for
$\langle\bar\rmR(X,Y)N,Z\rangle$.
\end{proof}

Using the fact that the vector field $\frac{\partial}{\partial t}$ is
parallel, we obtain the following equations.

\begin{prop}
For $X\in\mathfrak{X}(\cV)$ we have
$$\nabla_XT=\nu\rmS X,\quad
\rmd\nu(X)=-\langle\rmS X,T\rangle.$$
\end{prop}

\begin{proof}
We have $\frac{\partial}{\partial t}=T+\nu N$ and
$\bar\nabla_X\frac{\partial}{\partial t}=0$. Thus we get
$$0=\bar\nabla_XT+(\rmd\nu(X))N+\nu\bar\nabla_XN
=\nabla_XT+\langle\rmS X,T\rangle N+(\rmd\nu(X))N-\nu\rmS X.$$
Taking the tangential and the normal components in this equality, we
obtain the expected formulas.
\end{proof}

\begin{rem}
In the case of an orthonormal pair $(X,Y)$ we get
$$\langle\bar\rmR(X,Y)X,Y\rangle=
\kappa(1
-\langle Y,T\rangle^2-\langle X,T\rangle^2).$$
\end{rem}

The reader can also refer to section 3.2 in \cite{abresch}.

\subsection{Moving frames} \label{movingframes}

In this section we introduce some material about the technique of
moving frames. The reader can also refer to \cite{rosenberg}.

Let $\cV$ be a Riemannian manifold of dimension $n$, 
$\nabla$ its Levi-Civita connection, and
$\rmR$ the Riemannian curvature tensor. 
Let $\rmS$ be a field of symmetric operators
$\rmS_y:\rmT_y\cV\to\rmT_y\cV$.
Let $(e_1,\dots,e_n)$ be a local orthonormal frame on $\cV$ and
$(\omega^1,\dots,\omega^n)$ the dual basis of  $(e_1,\dots,e_n)$,
i.e., $$\omega^i(e_k)=\delta^i_k.$$ We also set
$$\omega^{n+1}=0.$$ 

We define the forms $\omega^i_j$, $\omega^{n+1}_j$, $\omega^i_{n+1}$
and $\omega^{n+1}_{n+1}$ on $\cV$ by
$$\omega^i_j(e_k)=\langle\nabla_{e_k}e_j,e_i\rangle,\quad
\omega^{n+1}_j(e_k)=\langle\rmS e_k,e_j\rangle,$$
$$\omega^j_{n+1}=-\omega^{n+1}_j,\quad
\omega^{n+1}_{n+1}=0.$$
Then we have
$$\nabla_{e_k}e_j=\sum_i\omega^i_j(e_k)e_i,\quad
\rmS e_k=\sum_j\omega^{n+1}_j(e_k)e_j.$$

Finally we set $R^i_{klj}=\langle\rmR(e_k,e_l)e_j,e_i\rangle$.

\begin{prop} \label{differentiation}
We have the following formulas:
\begin{equation} \label{diffomega1}
\rmd\omega^i+\sum_p\omega^i_p\wedge\omega^p=0,
\end{equation}
\begin{equation} \label{diffomega2}
\sum_p\omega^{n+1}_p\wedge\omega^p=0,
\end{equation}
\begin{equation} \label{diffomega3}
\rmd\omega^i_j+\sum_p\omega^i_p\wedge\omega^p_j=
-\frac{1}{2}\sum_k\sum_lR^i_{klj}\omega^k\wedge\omega^l,
\end{equation}
\begin{equation} \label{diffomega4}
\rmd\omega^{n+1}_j+\sum_p\omega^{n+1}_p\wedge\omega^p_j=
\frac{1}{2}\sum_k\sum_l\langle\nabla_{e_k}\rmS e_l
-\nabla_{e_l}\rmS e_k-\rmS[e_k,e_l],e_j\rangle\omega^k\wedge\omega^l.
\end{equation}
\end{prop}

\begin{proof}
These are well known formulas. However, since our conventions slightly
differ from those of \cite{tenenblat} and \cite{rosenberg}, we give a
proof for sake of clarity.

We have $\rmd\omega^i(e_p,e_q)=-\omega^i([e_p,e_q])
=-\omega^i(\nabla_{e_p}e_q-\nabla_{e_q}e_p)
=-\omega^i_q(e_p)+\omega^i_p(e_q)$ and
$\sum_k\omega^i_k\wedge\omega^k(e_p,e_q)
=\omega^i_q(e_p)-\omega^i_p(e_q)$, so \eqref{diffomega1} is proved.
And we have
$\sum_k(\omega^{n+1}_k\wedge\omega^k)(e_p,e_q)
=\omega^{n+1}_q(e_p)-\omega^{n+1}_p(e_q)
=\langle \rmS e_p,e_q\rangle-\langle \rmS e_q,e_p\rangle=0$, so
\eqref{diffomega2} is proved.

We have $\omega_j^i=\sum_k\langle e_i,\nabla_{e_k}e_j\rangle
\omega^k$, so
\begin{eqnarray*}
\rmd\omega_j^i & = & \sum_k\sum_l
e_l\langle e_i,\nabla_{e_k}e_j\rangle\omega^l\wedge\omega^k
+\sum_k\langle e_i,\nabla_{e_k}e_j\rangle
\rmd\omega^k \\
& = &
\sum_k\sum_l
(\langle\nabla_{e_l}e_i,\nabla_{e_k}e_j\rangle
+\langle e_i,\nabla_{e_l}\nabla_{e_k}e_j\rangle)
\omega^l\wedge\omega^k \\
& & -\sum_k\sum_l\langle e_i,\nabla_{e_k}e_j\rangle
\omega_l^k\wedge\omega^l.
\end{eqnarray*}
Moreover we have
\begin{eqnarray*}
\sum_k\sum_l\langle e_i,\nabla_{e_k}e_j\rangle
\omega_l^k\wedge\omega^l & = &
\sum_k\sum_l\sum_q
\langle e_i,\nabla_{e_k}e_j\rangle
\langle e_k,\nabla_{e_q}e_l\rangle
\omega^q\wedge\omega^l \\
& = &
\sum_l\sum_q
\langle e_i,\nabla_{\nabla_{e_q}e_l}e_j\rangle
\omega^q\wedge\omega^l.
\end{eqnarray*}
On the other hand we have
\begin{eqnarray*}
\sum_p\omega_p^i\wedge\omega_j^p & = &
\sum_k\sum_l\sum_p
\langle e_i,\nabla_{e_l}e_p\rangle
\langle e_p,\nabla_{e_k}e_j\rangle
\omega^l\wedge\omega^k \\
& = &
-\sum_k\sum_l\sum_p
\langle\nabla_{e_l}e_i,e_p\rangle
\langle e_p,\nabla_{e_k}e_j\rangle
\omega^l\wedge\omega^k \\
& = &
-\sum_k\sum_l
\langle\nabla_{e_l}e_i,\nabla_{e_k}e_j\rangle
\omega^l\wedge\omega^k.
\end{eqnarray*}
Thus we conclude that
$$\rmd\omega_j^i+\sum_p\omega_p^i\wedge\omega_j^p=
\sum_k\sum_l
\langle e_i,\nabla_{e_l}\nabla_{e_k}e_j
-\nabla_{\nabla_{e_l}e_k}e_j\rangle
\omega^l\wedge\omega^k.$$
Adding this equality with itself after exchanging $k$ and $l$ and
using the fact that $\omega^k\wedge\omega^l=
-\omega^l\wedge\omega^k$,  we get
\begin{eqnarray*}
2\left(\rmd\omega_j^i+\sum_p\omega_p^i\wedge\omega_j^p\right)
& = &
\sum_k\sum_l
\langle e_i,\rmR(e_k,e_l)e_j\rangle
\omega^l\wedge\omega^k,
\end{eqnarray*}
and finally we get \eqref{diffomega3}.

We have $\omega^{n+1}_j=\sum_k\langle\rmS e_k,e_j\rangle\omega^k$, so
\begin{eqnarray*}
\rmd\omega^{n+1}_j & = &
\sum_k\sum_l e_l\langle\rmS e_k,e_j\rangle
\omega^l\wedge\omega^k
+\sum_k\langle\rmS e_k,e_j\rangle\rmd\omega^k \\
& = &
\sum_k\sum_l(\langle\nabla_{e_l}\rmS e_k,e_j\rangle
+\langle\rmS e_k,\nabla_{e_l}e_j\rangle)\omega^l\wedge\omega^k
-\sum_k\sum_l\langle\rmS e_k,e_j\rangle
\omega_l^k\wedge\omega^l.
\end{eqnarray*}
Moreover we have
\begin{eqnarray*}
\sum_k\sum_l\langle\rmS e_k,e_j\rangle
\omega_l^k\wedge\omega^l & = &
\sum_k\sum_l\sum_q
\langle\rmS e_k,e_j\rangle
\langle e_k,\nabla_{e_q}e_l\rangle
\omega^q\wedge\omega^l \\
& = &
\sum_l\sum_q
\langle\rmS e_j,\nabla_{e_q}e_l\rangle
\omega^q\wedge\omega^l.
\end{eqnarray*}
On the other hand we have
\begin{eqnarray*}
\sum_p\omega^{n+1}_p\wedge\omega_j^p & = &
\sum_k\sum_p\langle\rmS e_k,e_p\rangle
\omega^k\wedge\omega_j^p \\
& = &
\sum_k\sum_p\sum_l
\langle\rmS e_k,e_p\rangle
\langle e_p,\nabla_{e_l}e_j\rangle\omega^k\wedge\omega^l \\
& = &
\sum_k\sum_l
\langle\rmS e_k,\nabla_{e_l}e_j\rangle\omega^k\wedge\omega^l.
\end{eqnarray*}
Thus we conclude that
\begin{eqnarray*}
\rmd\omega^{n+1}_j+\sum_p\omega^{n+1}_p\wedge\omega_j^p & = &
\sum_k\sum_l(\langle\nabla_{e_l}\rmS e_k,e_j\rangle
-\langle\rmS e_j,\nabla_{e_l}e_k\rangle)
\omega^l\wedge\omega^k \\
& = &
\sum_k\sum_l\langle e_j,\nabla_{e_l}\rmS e_k
-\rmS\nabla_{e_l}e_k\rangle
\omega^l\wedge\omega^k. \\
\end{eqnarray*}
Adding this equality with itself after exchanging $k$ and $l$ and
using the fact that $\omega^k\wedge\omega^l=
-\omega^l\wedge\omega^k$, we get
\begin{eqnarray*}
2\left(\rmd\omega^{n+1}_j+\sum_p\omega^{n+1}_p\wedge\omega_j^p\right)
& = &
\sum_k\sum_l\langle e_j,\nabla_{e_l}\rmS e_k
-\nabla_{e_k}\rmS e_l-\rmS[e_l,e_k]\rangle
\omega^l\wedge\omega^k,
\end{eqnarray*}
and finally we get \eqref{diffomega4}.
\end{proof}

\subsection{Some facts about hypersurfaces of $\s^n\times\R$ and
$\h^n\times\R$} \label{hypersurfaces}

In this section we consider an orientable hypersurface $\cV$ of
$\M^n\times\R$ with $\M^n=\s^n$ or $\M^n=\h^n$.

We denote by $\LL^p$ the $p$-dimensional Lorentz space, i.e., $\R^p$
endowed with the quadradic form 
$$-(\rmd x^0)^2+(\rmd x^1)^2+\dots+(\rmd x^{p-1})^2.$$

We will use the following inclusions:
we have 
$$\s^n=\{(x^0,\dots,x^n)\in\R^{n+1};(x^0)^2+\sum_i(x^i)^2=1\}$$
and so
$$\s^n\times\R\subset\R^{n+1}\times\R=\R^{n+2},$$
and we have 
$$\h^n=\{(x^0,\dots,x^n)\in\LL^{n+1};-(x^0)^2+\sum_i(x^i)^2=-1,x^0>0\}$$
and so
$$\h^n\times\R\subset\LL^{n+1}\times\R=\LL^{n+2}.$$
In the case of $\s^n\times\R$ we set $\kappa=1$ and
$\E^{n+2}=\R^{n+2}$. In the case of $\h^n\times\R$ we set 
$\kappa=-1$ and $\E^{n+2}=\LL^{n+2}$.

We denote by $\nabla$, $\bar\nabla$ and $\bar{\bar\nabla}$ the
connections of $\cV$, $\M^n\times\R$ and $\E^{n+2}$ respectively, by
$\bar N(x)$ the 
normal to $\M^n\times\R$ in $\E^{n+2}$ at a point $x\in\M^n\times\R$, i.e., 
$$\bar N(x)=(x^0,\dots,x^n,0),$$ 
and by $N(x)$ the normal to $\cV$ in
$\M^n\times\R$ at a point $x\in\cV$.
We denote by $\rmS$ the shape operator of $\cV$ in $\M^n\times\R$.
The shape operator of $\M^n\times\R$ is
$\bar\rmS X=-\kappa\rmd\bar N(X)=\kappa\left(-X+
\left\langle X,\frac{\partial}{\partial t}\right\rangle
\frac{\partial}{\partial t}\right)$; we should be careful with the
sign convention in the definition of the shape operator: here we have
chosen
$$\bar{\bar\nabla}_XY=\bar\nabla_XY+\langle\bar\rmS X,Y\rangle\bar N,$$
i.e., 
$$\langle\bar\rmS X,Y\rangle
=\kappa\langle\bar{\bar\nabla}_XY,\bar N\rangle,$$
because in the case of $\s^n\times\R$ we have $\langle\bar
N,\bar N\rangle=1$ whereas in the case of $\h^n\times\R$ we have
$\langle\bar N,\bar N\rangle=-1$.

Let $(e_1,\dots,e_n)$ be a local
orthonormal frame on $\cV$, $e_{n+1}=N$ and $e_0=\bar N$ (on $\cV$).
We define the forms $\omega^i_j$, $\omega^{n+1}_j$,
$\omega^i_{n+1}$ and $\omega^{n+1}_{n+1}$ as in section
\ref{movingframes}. Moreover we set 
$$\omega^0_\gamma(e_k)=\langle\bar\rmS e_k,e_\gamma\rangle
=-\kappa\langle e_k,e_\gamma\rangle+\kappa
\left\langle e_k,\frac{\partial}{\partial t}\right\rangle
\left\langle e_\gamma,\frac{\partial}{\partial t}\right\rangle,$$
$$\omega^\gamma_0=-\kappa\omega^0_\gamma.$$
With these definitions we have
$$\bar{\bar\nabla}_{e_k}e_\beta
=\sum_\alpha\omega^\alpha_\beta(e_k)e_\alpha.$$

Let $(E_0,\dots,E_{n+1})$ be the canonical frame of
$\E^{n+2}$ (with $\langle E_0,E_0\rangle=\kappa$ and
$E_{n+1}=\frac{\partial}{\partial t}$). Let $A\in\cM_{n+2}(\R)$ be the
matrix (the indices going from $0$ to $n+1$) whose columns are the 
coordinates of the $e_\beta$ in the frame $(E_\alpha)$, i.e.,
$$e_\beta=\sum_\alpha A^\alpha_\beta E_\alpha.$$
Then, on the one hand we have
$$\bar{\bar\nabla}_{e_k}e_\beta
=\sum_\alpha\rmd A^\alpha_\beta(e_k)E_\alpha,$$
and on the other hand we have
$$\bar{\bar\nabla}_{e_k}e_\beta=
\sum_\alpha\sum_\gamma\omega^\gamma_\beta(e_k)
A^\alpha_\gamma E_\alpha.$$ 
Thus we have
$$A^{-1}\rmd A=\Omega$$ with
$\Omega=(\omega^\alpha_\beta)\in\cM_{n+2}(\R)$, the indices going from
$0$ to $n+1$.

Setting $G=\diag(\kappa,1,\dots,1)\in\cM_{n+2}(\R)$, we have
$$A\in\mathrm{SO}^+(\E^{n+2}),\quad\Omega\in\mathfrak{so}(\E^{n+2})$$
where $\mathrm{SO}^+(\E^{n+2})$ is the connected component of
$\mathrm{I}_{n+2}$ in
$$\mathrm{SO}(\E^{n+2})=\{Z\in\cM_{n+2}(\R);\tra ZGZ=G,\det Z=1\},$$
and where $$\mathfrak{so}(\E^{n+2})=
\{H\in\cM_{n+2}(\R);\tra HG+GH=0\}.$$
In the case of $\s^n\times\R$ we have
$\mathrm{SO}^+(\E^{n+2})=\mathrm{SO}(\R^{n+2})$.

\section{Isometric immersions into $\s^n\times\R$ and $\h^n\times\R$}

\subsection{The compatibility equations}

We consider a simply connected Riemannian manifold $\cV$ of dimension
$n$. Let $\rmd s^2$ be the metric on $\cV$ (we will also denote it by
$\langle\cdot,\cdot\rangle$), $\nabla$  
the Riemannian connection of $\cV$ and $\rmR$ its Riemann curvature
tensor. Let $\rmS$ be a field of symmetric operators
$\rmS_y:\rmT_y\cV\to\rmT_y\cV$, $T$ a vector field on $\cV$
such that $||T||\leqslant 1$ and $\nu$ a smooth function on $\cV$ 
such that $\nu^2\leqslant 1$.

The compatibility equations for hypersurfaces in $\s^n\times\R$ and
$\h^n\times\R$ 
established in section \ref{compatibilitymtimesr} suggest to introduce
the following definition. 

\begin{defn}
We say that $(\rmd s^2,\rmS,T,\nu)$ satisfies the
compatibility equations respectively for $\s^n\times\R$ and
$\h^n\times\R$
if $$||T||^2+\nu^2=1$$ and, for all $X,Y,Z\in\mathfrak{X}(\cV)$, 
\begin{equation} \label{gaussmtimesr}
\begin{array}{lll}
\rmR(X,Y)Z & = &
\langle\rmS X,Z\rangle\rmS Y
-\langle\rmS Y,Z\rangle\rmS X \\
& & +\kappa(\langle X,Z\rangle Y-
\langle Y,Z\rangle X
-\langle Y,T\rangle\langle X,Z\rangle T  \\
& & -\langle X,T\rangle\langle Z,T\rangle Y
+\langle X,T\rangle\langle Y,Z\rangle T
+\langle Y,T\rangle\langle Z,T\rangle X),
\end{array}
\end{equation}
\begin{equation} \label{codazzimtimesr}
\nabla_X\rmS Y-\nabla_Y\rmS X-\rmS[X,Y]=
\kappa\nu(\langle Y,T\rangle X-\langle X,T\rangle Y),
\end{equation}
\begin{equation} \label{conditionT1}
\nabla_XT=\nu\rmS X,
\end{equation}
\begin{equation} \label{conditionT2}
\rmd\nu(X)=-\langle\rmS X,T\rangle,
\end{equation}
where $\kappa=1$ and $\kappa=-1$ for $\s^n\times\R$ and $\h^n\times\R$
respectively.
\end{defn}

\begin{rem}
We notice that \eqref{conditionT1} implies \eqref{conditionT2} except
when $\nu=0$ (by differentiating the identity $\langle T,T\rangle
+\nu^2=1$ with respect to $X$).
\end{rem}

%
%
%
%
%
%
%
%

\subsection{Codimension $1$ isometric immersions into $\s^n\times\R$
and $\h^n\times\R$} \label{sectiontheorem}

In this section we will prove the following theorem.

\begin{thm} \label{isometry}
Let $\cV$ be a simply connected Riemannian manifold of dimension $n$,
$\rmd s^2$ its metric and $\nabla$ its 
Riemannian connection. Let $\rmS$ be a field of symmetric operators
$\rmS_y:\rmT_y\cV\to\rmT_y\cV$, $T$ a vector field on $\cV$ and $\nu$
a smooth function on $\cV$ such that $||T||^2+\nu^2=1$. 

Let $\M^n=\s^n$ or $\M^n=\h^n$.
Assume that $(\rmd s^2,\rmS,T,\nu)$ satisfies the
compatibility equations for $\M^n\times\R$. Then
there exists an isometric immersion $f:\cV\to\M^n\times\R$
such that
the shape operator with respect to the normal $N$ associated to $f$ is
$$\rmd f\circ\rmS\circ\rmd f^{-1}$$ and such that
$$\frac{\partial}{\partial t}=\rmd f(T)+\nu N.$$ Moreover the
immersion is unique up to a global isometry of $\M^n\times\R$
preserving the orientations of both $\M^n$ and $\R$.
\end{thm}

To prove this theorem, we consider a local orthonormal frame
$(e_1,\dots,e_n)$ on $\cV$ and the forms $\omega^i$, $\omega^{n+1}$,
$\omega^i_j$, $\omega^{n+1}_j$, $\omega^i_{n+1}$ and
$\omega^{n+1}_{n+1}$ as in section \ref{movingframes}.
We set $\E^{n+2}=\R^{n+2}$ or $\E^{n+2}=\LL^{n+2}$ (according to
$\M^n$). We denote by 
$(E_0,\dots,E_{n+1})$ the canonical frame of $\E^{n+2}$ (with $\langle
E_0,E_o\rangle=-1$ in the case of $\LL^{n+2}$); in particular we
have $E_{n+1}=\frac{\partial}{\partial t}$. We set 
$$T^k=\langle T,e_k\rangle,\quad T^{n+1}=\nu,\quad
T^0=0.$$

Moreover we set
$$\omega^0_j(e_k)=\kappa(T^jT^k-\delta^k_j),\quad
\omega^0_{n+1}(e_k)=\kappa\nu T^k,$$
$$\omega^i_0=-\kappa\omega^0_i,\quad
\omega^{n+1}_0=-\kappa\omega^0_{n+1},\quad
\omega^0_0=0.$$ 
We define the one-form $\eta$ on $\cV$ by
$$\eta(X)=\langle T,X\rangle.$$
In the frame $(e_1,\dots,e_n)$ we have $\eta=\sum_kT^k\omega^k$.
Finaly we define the following matrix of one-forms:
$$\Omega=(\omega^\alpha_\beta)\in\cM_{n+2}(\R),$$
the indices going from $0$ to $n+1$.

From now on we assume that the hypotheses of theorem \ref{isometry}
are satisfied. We first prove some technical lemmas that are
consequences of the compatibility equations.

\begin{lemma} \label{diffeta}
We have $$\rmd\eta=0.$$
\end{lemma}

\begin{proof}
We have
\begin{eqnarray*}
\rmd\eta(X,Y) & = & X\cdot\eta(Y)-Y\cdot\eta(X)-\eta([X,Y]) \\
& = & \langle\nabla_XT,Y\rangle-\langle\nabla_YT,X\rangle \\
& = & \langle\nu\rmS X,Y\rangle-\langle\nu\rmS Y,X\rangle \\
& = & 0.
\end{eqnarray*}
We have used condition \eqref{conditionT1}.
\end{proof}

\begin{lemma} \label{diffT}
We have 
$$\rmd T^\alpha=\sum_\gamma T^\gamma\omega^\gamma_\alpha.$$
\end{lemma}

\begin{proof}
This is a consequence of condition \eqref{conditionT1}
for $\alpha=j$, of condition \eqref{conditionT2} for $\alpha=n+1$, and
of the definitions for $\alpha=0$.
\end{proof}

\begin{lemma} \label{diffOmega}
We have $$\rmd\Omega+\Omega\wedge\Omega=0.$$
\end{lemma}

\begin{proof}
We set $\Psi=\rmd\Omega+\Omega\wedge\Omega$ and 
$R^i_{klj}=\langle\rmR(e_k,e_l)e_j,e_i\rangle$.

By proposition \ref{differentiation} we have
$$\Psi^i_j=-\frac{1}{2}\sum_k\sum_lR^i_{klj}\omega^k\wedge\omega^l
+\omega^i_{n+1}\wedge\omega^{n+1}_j
+\omega^i_0\wedge\omega^0_j.$$
Since the Gauss equation \eqref{gaussmtimesr} is satisfied we have
$$R^i_{klj}
=\bar R^i_{klj}+\omega^{n+1}_j\wedge\omega^{n+1}_i(e_k,e_l)$$
with $$\bar R^i_{klj}
=\kappa(\delta^k_j\delta^l_i
-\delta^l_j\delta^k_i-T^lT^i\delta^k_j
-T^kT^j\delta^l_i+T^kT^i\delta^l_j+T^lT^j\delta^k_i).$$ On the other
hand, a computation shows that 
$\omega^i_0\wedge\omega^0_j(e_k,e_l)=\bar R^i_{klj}$. Thus we
have $R^i_{klj}=\omega^i_{n+1}\wedge\omega^{n+1}_j(e_k,e_l)
+\omega^i_0\wedge\omega^0_j(e_k,e_l)$. We conclude that
$\Psi^i_j=0$.

By proposition \ref{differentiation} we have
$$\Psi^{n+1}_j=\frac{1}{2}\sum_k\sum_l
\langle\nabla_{e_k}\rmS e_l-\nabla_{e_l}\rmS e_k-\rmS[e_k,e_l],
e_j\rangle\omega^k\wedge\omega^l+
\omega^{n+1}_0\wedge\omega^0_j.$$
Since the Codazzi equation \eqref{codazzimtimesr} is satisfied we have
$$\langle\nabla_{e_k}\rmS e_l-\nabla_{e_l}\rmS e_k-\rmS[e_k,e_l],
e_j\rangle
=\kappa(T^lT^{n+1}\delta^k_j-T^kT^{n+1}\delta^l_j).$$ On the
other hand, a computation shows that
$\omega^{n+1}_0\wedge\omega^0_j(e_k,e_l)=
\kappa(T^kT^{n+1}\delta^l_j-T^lT^{n+1}\delta^k_j)$. We conclude that
$\Psi^{n+1}_j=0$.

We have $\omega^0_j=\kappa(T^j\eta-\omega^j)$. Since $\rmd\eta=0$ (by
lemma \ref{diffeta}) we get
$$\rmd\omega^0_j=\kappa(\rmd T^j\wedge\eta-\rmd\omega^j)
=\kappa\rmd T^j\wedge\eta+\kappa\sum_k\omega^j_k\wedge\omega^k$$ 
by proposition
\ref{differentiation}. Thus by a straightforward computation we get
\begin{eqnarray*}
\Psi^0_j(e_p,e_q) & = & \rmd\omega^0_j(e_p,e_q)
+\sum_k\omega^0_k\wedge\omega^k_j(e_p,e_q)
+\omega^0_{n+1}\wedge\omega^{n+1}_j(e_p,e_q) \\
& = & \kappa(\rmd T^j(e_p)\eta(e_q)-\rmd T^j(e_q)\eta(e_p)
+\omega^j_q(e_p)-\omega^j_p(e_q)) \\
& & +\kappa\left(T^p\sum_kT^k\omega^k_j(e_q)
-T^q\sum_kT^k\omega^k_j(e_p)
-\omega^p_j(e_q)+\omega^q_j(e_p)\right) \\
& & +\kappa\left(T^pT^{n+1}\omega^{n+1}_j(e_q)
-T^qT^{n+1}\omega^{n+1}_j(e_p)\right).
\end{eqnarray*}
Using the definition of $\eta$ and lemma \ref{diffT} for $\alpha=j$,we
conclude that $\hat\Psi^{n+2}_j=0$. 

We have $\omega^0_{n+1}=\kappa T^{n+1}\eta$, and so
$\rmd\omega^0_{n+1}=\kappa\rmd T^{n+1}\wedge\eta$ by lemma
\ref{diffeta}. Thus by a straightforward computation we get
\begin{eqnarray*}
\Psi^0_{n+1}(e_p,e_q) & = & \rmd\omega^0_{n+1}(e_p,e_q)
+\sum_k\omega^0_k\wedge\omega^k_{n+1}(e_p,e_q) \\
& = & \kappa(T^q\rmd T^{n+1}(e_p)-T^p\rmd T^{n+1}(e_q)) \\
& & +\kappa\left(T^p\sum_kT^k\omega^k_{n+1}(e_q)
-T^q\sum_kT^k\omega^k_{n+1}(e_p)\right) \\
& & +\kappa(-\omega^p_{n+1}(e_q)+\omega^q_{n+1}(e_p)).
\end{eqnarray*}
The last two terms cancel because $\rmS$ is symmetric. Using lemma
\ref{diffT} for $\alpha=n+1$, we conclude that $\Psi^0_{n+1}=0$.

The fact that $\Psi^0_0=0$ and $\Psi^{n+1}_{n+1}=0$ is clear. We
conclude by noticing that $\Psi^i_{n+1}=-\Psi^{n+1}_i=0$.
\end{proof}

%

For $y\in\cV$, let $\cZ(y)$ be the set of matrices
$Z\in\mathrm{SO}^+(\E^{n+2})$ such 
that the coefficients of the last line of $Z$ are the $T^\beta(y)$.
It is a manifold of dimension $\frac{n(n+1)}{2}$ (since the map 
$F:\mathrm{SO}^+(\E^{n+2})\to\s(\E^{n+2}),Z\mapsto(Z^{n+1}_\beta)_\beta$
(i.e., $F(Z)$ is the last line of $Z$), where
$\s(\E^{n+2})=\{x\in\E^{n+2};\langle E,E\rangle=1\}$, is a submersion).

We now prove the following proposition.

\begin{prop} \label{matrixA}
Assume that the compatibility equations for $\M^n\times\R$ are
satisfied. Let $y_0\in\cV$ and $A_0\in\cZ(y_0)$. Then there exist a
neighbourhood $U_1$ of
$y_0$ in $\cV$ and a unique map $A:U_1\to\mathrm{SO}^+(\E^{n+2})$ such that
$$A^{-1}\rmd A=\Omega,$$
$$\forall y\in U_1,\quad A(y)\in\cZ(y),$$
$$A(y_0)=A_0.$$
\end{prop}

\begin{proof}
Let $U$ be a coordinate neighbourhood in $\cV$. The set 
$$\cF=\{(y,Z)\in U\times\mathrm{SO}^+(\E^{n+2});Z\in\cZ(y)\}$$ is a
manifold of dimension $n+\frac{n(n+1)}{2}$, and
$$\rmT_{(y,Z)}\cF=\{(u,\zeta)\in\rmT_yU\oplus\rmT_Z\mathrm{SO}^+(\E^{n+2});
\zeta^{n+1}_\beta=(\rmd T^\beta)_y(u)\}.$$
Indeed, in the neighbourhood of point of $U$ there exists a map
$y\mapsto M(y)\in\mathrm{SO}^+(\E^{n+2})$ such that the last line of
$M(y)$ is $(T^\beta(y))_\beta$, and we have $Z\in\cZ(y)$ if and only
if 
$$ZM(y)^{-1}=\left(\begin{array}{cc}
B & 0 \\
0 & 1
\end{array}\right)$$ for some $B\in\mathrm{SO}^+(\E^{n+1})$; then, if
$\varphi$ is a local parametrization of the set of such matrices, the
map $(y,v)\mapsto(y,\varphi(v)M(y))$ is a local parametrization of
$\cF$.

Let $Z$ denote the projection
$U\times\mathrm{SO}^+(\E^{n+2})
\to\mathrm{SO}^+(\E^{n+2})\subset\cM_{n+2}(\R)$.
We consider on $\cF$ the 
following matrix of $1$-forms:
$$\Theta=Z^{-1}\rmd Z-\Omega,$$ namely for $(y,Z)\in\cF$ we have
$$\Theta_{(y,Z)}:\rmT_{(y,Z)}\cF\to\cM_{n+2}(\R),$$
$$\Theta_{(y,Z)}(u,\zeta)=Z^{-1}\zeta-\Omega_y(u).$$

We claim that, for each $(y,Z)\in\cF$, the space
$$\cD(y,Z)=\ker\Theta_{(y,Z)}$$ has dimension $n$.
We first notice that
the matrix $\Theta$ belongs to $\mathfrak{so}(\E^{n+2})$ since
$\Omega$ and $Z^{-1}\rmd Z$ do. Moreover we have
$$(Z\Theta)^{n+1}_\beta
=\rmd Z^{n+1}_\beta-\sum_\gamma Z^{n+1}_\gamma\omega^\gamma_\beta
=\rmd T^\beta-\sum_\gamma T^\gamma\omega^\gamma_\beta
=0$$ by lemma \ref{diffT}. Thus the values of
$\Theta_{(y,Z)}$ lie in the space
$$\cH=\{H\in\mathfrak{so}^+(\E^{n+2});(ZH)^{n+1}_\beta=0\},$$
which has dimension $\frac{n(n+1)}{2}$ (indeed, the map 
$F:\mathrm{SO}^+(\E^{n+2})\to\s(\E^{n+2}),
Z\mapsto(Z^{n+1}_\beta)_\beta$ is a submersion, and we
have $H\in\cH$ if and 
only if $ZH\in\ker(\rmd F)_Z$).
Moreover, the space $\rmT_{(y,Z)}\cF$ contains the subspace
$\{(0,ZH);H\in\cH\}$, and the restriction of $\Theta_{(y,Z)}$ on this
subspace is the map $(0,ZH)\mapsto H$. Thus
$\Theta_{(y,Z)}$ is onto 
$\cH$, and consequently the linear map $\Theta_{(y,Z)}$ has rank
$\frac{n(n+1)}{2}$. This finishes proving the claim.

We now prove that the distribution $\cD$ is involutive.
Using lemma \ref{diffOmega} we get
\begin{eqnarray*}
\rmd\Theta & = & -Z^{-1}\rmd Z\wedge Z^{-1}\rmd Z-\rmd\Omega \\
& = & -(\Theta+\Omega)\wedge(\Theta+\Omega)-\rmd\Omega \\
& = & -\Theta\wedge\Theta-\Theta\wedge\Omega-\Omega\wedge\Theta.
\end{eqnarray*}
From this formula we deduce that if $\xi_1,\xi_2\in\cD$, then
$\rmd\Theta(\xi_1,\xi_2)=0$, and so
$\Theta([\xi_1,\xi_2])=\xi_1\cdot\Theta(\xi_2)
-\xi_2\cdot\Theta(\xi_1)-\rmd\Theta(\xi_1,\xi_2)=0$, i.e.,
$[\xi_1,\xi_2]\in\cD$. Thus the distribution $\cD$ is involutive, and
so, by the theorem of Frobenius, it is integrable.

Let $\cA$ be the integral manifold through $(y_0,A_0)$.
If $\zeta\in\rmT_{A_0}\mathrm{SO}^+(\E^{n+2})$ is such that
$(0,\zeta)\in\rmT_{(y_0,A_0)}\cA=\cD(y_0,A_0)$, 
then we have $0=\Theta_{(y_0,A_0)}(0,\zeta)=A_0^{-1}\zeta$. This
proves that
$$\rmT_{(y_0,A_0)}\cA\cap
\left(\{0\}\times\rmT_{A_0}\mathrm{SO}^+(\E^{n+2})\right)=\{0\}.$$
Thus the manifold
$\cA$ is locally the graph of a function
$A:U_1\to\mathrm{SO}^+(\E^{n+2})$ where $U_1$ is a neighbourhood of
$y_0$ in $U$. By construction, this map satisfies the properties of
proposition \ref{matrixA} and is unique. 
\end{proof}

We now prove the theorem.

\begin{proof}[Proof of theorem \ref{isometry}]
Let $y_0\in\cV$, $A\in\cZ(y_0)$ and $t_0\in\R$.
We consider on $\cV$ a local orthonormal frame $(e_1,\dots,e_n)$ in
the neighbourhood of $y_0$ and we keep the same notations. Then by
proposition \ref{matrixA} there exists a unique map 
$A:U_1\to\mathrm{SO}^+(\E^{n+2})$ such that
$$A^{-1}\rmd A=\Omega,$$
$$\forall y\in U_1,\quad A(y)\in\cZ(y),$$
$$A(y_0)=A_0,$$ where $U_1$ is a neighbourhood of $y_0$, which we can
assume simply connected.

We set $f^0=A^0_0$, $f^i=A^i_0$ and we call $f^{n+1}$
the unique function on $U_1$ such that $\rmd f^{n+1}=\eta$ and
$f^{n+1}(y_0)=t_0$ (this function exists since $U_1$ is simply
connected and $\rmd\eta=0$). Thus we defined a map $f:U_1\to\E^{n+2}$. Since
$A^{n+1}_0=T^0=0$ and $A\in\mathrm{SO}^+(\E^{n+2})$, in the case of
$\s^n\times\R$ we have
$(f^0)^2+\sum_i(f^i)^2=\sum_\alpha(A^\alpha_0)^2=1$, and in the case
of $\h^n\times\R$ we have $-(f^0)^2+\sum_i(f^i)^2=
-(A^0_0)^2+\sum_i(A^i_0)^2+(A^{n+1}_0)^2=-1$ and $f^0=A^0_0>0$. Thus
in both cases we have $(f^0,\dots,f^n)\in\M^n$, i.e., 
the values of $f$ lie in $\M^n\times\R$.

Since $\rmd A=A\Omega$, we have, for $\alpha<n+1$, 
\begin{eqnarray*}
\rmd f^\alpha(e_k) & = & \sum_jA^\alpha_j\omega^j_0(e_k)
+A^\alpha_{n+1}\omega^{n+1}_0(e_k) \\
& = & \sum_jA^\alpha_j(\delta^k_j-T^jT^k)-A^\alpha_{n+1}T^{n+1}T^k \\
& = & A^\alpha_k-T^k\sum_\beta A^\alpha_\beta A^{n+1}_\beta \\
& = & A^\alpha_k,
\end{eqnarray*}
and $$\rmd f^{n+1}(e_k)=\eta(e_k)=T^k=A^{n+1}_k.$$ This means that
$\rmd f(e_k)$ is given by the colunm $k$ of the matrix $A$.

Since $A$ is an invertible matrix, $\rmd f$ has rank $n$, and
so $f$ is an immersion. And since 
$A\in\mathrm{SO}^+(\E^{n+2})$, we have 
$\langle\rmd f(e_p),\rmd f(e_q)\rangle
=\delta^p_q$, and so $f$ is an isometry.

The columns of $A(y)$ form a direct orthonormal frame of
$\E^{n+2}$. Columns $1$ to $n$ form a direct orthonormal frame of
$\rmT_{f(y)}f(\cV)$ and column $0$
is the projection of $f(y)$ on $\M^n\times\{0\}$, i.e.,
the unit normal $\bar N(f(y))$ to $\M^n\times\R$ at the point
$f(y)$. Thus column $(n+1)$ is the unit normal $N(f(y))$ to
$f(\cV)$ in $\M^n\times\R$ at the point $f(y)$. 

We set $X_j=\rmd f(e_j)$. Then we have
\begin{eqnarray*}
\langle\rmd X_j(X_k),N\rangle & = &
\sum_\alpha\rmd A^\alpha_j(e_k)A^\alpha_{n+1}
=\sum_\alpha\sum_\gamma A^\alpha_\gamma
A^\alpha_{n+1}\omega^\gamma_j(e_k) \\
& = & \omega^{n+1}_j(e_k)=\langle\rmS e_k,e_j\rangle.
\end{eqnarray*}
This means that the shape operator of $f(\cV)$ in $\M^n\times\R$ is 
$\rmd f\circ\rmS\circ\rmd f^{-1}$.

Finally, the coefficients of the vertical vector
$\frac{\partial}{\partial t}=E_{n+1}$ in 
the orthonormal frame $(\bar N,X_1,\dots,X_n,N)$ are given by the last
line of $A$. Since $A(y)\in\cZ(y)$ for all $y\in U_2$ we get
$$\frac{\partial}{\partial t}=\sum_j T^jX_j+T^{n+1}N
=\rmd f(T)+\nu N.$$

We now prove that the local immersion is unique up to a global
isometry of $\M^n\times\R$. Let $\tilde f:U_3\to\M^n\times\R$ be
another immersion satisfying the conclusion of the theorem, where
$U_3$ is a simply connected neighbourhood of $y_0$ included in $U_1$,
let $(\tilde X_\beta)$ be the associated frame (i.e., $\tilde
X_j=\rmd\tilde f(e_j)$, $\tilde X_{n+1}$ is the normal of $\tilde
f(\cV)$ in $\M^n\times\R$ and $\tilde X_0$ is the normal to
$\M^n\times\R$ in $\E^{n+2}$) and let $\tilde A$ the
matrix of the coordinates of the frame $(\tilde X_\beta)$ in the
frame $(E_\alpha)$. Up to a direct isometry of $\M^n\times\R$,
we can
assume that $f(y_0)=\tilde f(y_0)$ and that the frames
$(X_\beta(y_0))$ and $(\tilde X_\beta(y_0))$ coincide, i.e.,
$A(y_0)=\tilde A(y_0)$. We notice that this isometry necessarily fixes
$\frac\partial{\partial t}$ since the $T^\alpha$ are the same for $x$
and $\tilde x$. The matrices $A$ and $\tilde A$ satisfy
$A^{-1}\rmd A=\Omega$ and $\tilde A^{-1}\rmd\tilde A=\Omega$ (see
section \ref{hypersurfaces}),
$A(y),\tilde A(y)\in\cZ(y)$ and $A(y_0)=\tilde A(y_0)$,
thus by the uniqueness of the solution of the equation in proposition
\ref{matrixA} we get $A(y)=\tilde A(y)$. Considering the columns $0$
of these matrices, we get $f^i=\tilde f^i$ and $f^0=\tilde
f^0$. Finally we have $\rmd f^{n+1}=\eta=\rmd\tilde f^{n+1}$ and
$f^{n+1}(y_0)=\tilde f^{n+1}(y_0)$, thus we have $f^{n+1}=\tilde
f^{n+1}$. This finishes proving that $f=\tilde f$ on $U_3$. 

Finally we prove that this local immersion $f$ can be extended to
$\cV$ in a unique way. Let $y_1\in\cV$. Then there exists a curve
$\Gamma:[0,1]\to\cV$ such that $\Gamma(0)=y_0$ and
$\Gamma(1)=y_1$. Each point of $\Gamma$ has a neighbourhood such that
there exists an isometric immersion (unique up to an isometry of
$\M^n\times\R$ preserving the orientations of $\M^n$ and $\R$) of this
neighbourhood satisfying the properties of the 
theorem. From this family of neighbourhoods we can extract a finite
family $(W_1,\dots,W_p)$ covering $\Gamma$ with $W_1=U_1$. Then the
above uniqueness argument shows that we can extend successively the
immersion $f$ to the $W_k$ in a unique way. In particular $f(y_1)$ is
defined. Moreover, this value $f(y_1)$ does not depend on the choice
of the curve $\Gamma$ joining $y_0$ to $y_1$ because $\cV$ is simply
connected.
\end{proof}

\begin{prop} \label{changeofsigns}
If $(\rmd s^2,\rmS,T,\nu)$ satisfies the compatibilty equations and
correspond to an immerion $f:\Sigma\to\M^n\times\R$, then 
$(\rmd s^2,-\rmS,T,-\nu)$, $(\rmd s^2,-\rmS,-T,\nu)$ and 
$(\rmd s^2,\rmS,-T,-\nu)$ also satisfy the compatibilty equations and
they correspond to the immersion $\sigma\circ f$ where $\sigma$ is an
isometry of $\M^n\times\R$
\begin{enumerate}
\item reversing the orientation of $\M^n$ and
preserving the orientation of $\R$ in the case of $(\rmd
s^2,-\rmS,T,-\nu)$,
\item preserving the
orientation of $\M^n$ and reversing the orientation of $\R$
in the case of $(\rmd s^2,-\rmS,-T,\nu)$,
\item reversing the orientations of both $\M^n$ and $\R$ in the case
of $(\rmd s^2,\rmS,-T,-\nu)$.
\end{enumerate}
\end{prop}

\begin{proof}
We deal with the first case (the two others are similar). Let $\hat
f=\sigma\circ f$. Then the normal to $\M^n\times\R$ is
$\sigma\circ\bar N$, and since $\sigma$ reverses the orientation of
$\M^n\times\R$ the normal to $\hat f(\cV)$ in $\M^n\times\R$ is
$\hat N=-\sigma\circ N$. From this we deduce that $\hat\rmS=-\rmS$. Moreover
we have $\frac{\partial}{\partial t}=\rmd f(T)+\nu N$, and so, since
$\sigma$ preserves the orientation of $\R$ we have 
$$\frac{\partial}{\partial t}=\sigma\circ\rmd f(T)+\nu\sigma\circ N
=\rmd\hat f(T)-\nu\hat N.$$ We conclude that $\hat T=T$ and
$\hat\nu=-\nu$.
\end{proof}


\subsection{Remark: another proof in the case of $\h^n\times\R$}

In this section we outline another proof of theorem \ref{isometry} in
the case of $\h^n\times\R$ that does not involve the Lorentz
space. Greek letters will denote indices between $1$ and $n+1$.

We first consider an orientable hypersurface $\cV$ of an
$(n+1)$-dimensionnal Riemannian manifold $\bar\cV$.
Let $(e_1,\dots,e_n)$ be a local 
orthonormal frame on $\cV$, $e_{n+1}$ the normal to $\cV$, and
$(E_1,\dots,E_{n+1})$ a local orthonormal frame on $\bar\cV$. We denote
by $\nabla$ and $\bar\nabla$ the Riemannian connections on
$\cV$ and $\bar\cV$ respectively, and by $\rmS$ the shape operator of
$\cV$ (with respect to the normal $e_{n+1}$). We define the forms
$\omega^\alpha$, $\omega^\alpha_\beta$ on $\cV$ as in
section \ref{movingframes}. Then we have
$$\bar\nabla_{e_k}e_\beta=\sum_\gamma\omega^\gamma_\beta(e_k)e_\gamma.$$

Let $A\in\mathrm{SO}_{n+1}(\R)$ be the matrix whose columns are the
coordinates of the $e_\beta$ in the frame $(E_\alpha)$, namely
$A^\alpha_\beta=\langle e_\beta,E_\alpha\rangle$. Let
$\Omega=(\omega^\alpha_\beta)\in\cM_{n+1}(\R)$. The matrix $A$
satisfies the following equation:
$$A^{-1}\rmd A=\Omega+L(A)$$
with $$L(A)^\alpha_\beta=\sum_k
\left(\sum_{\gamma,\delta,\varepsilon}A^\varepsilon_\alpha
A^\gamma_kA^\delta_\beta
\bar\Gamma_{\gamma\alpha}^\delta\right)\omega^k,$$
where the $\bar\Gamma_{\gamma\alpha}^\delta$ are the Christoffel
symbols of the frame $(E_\alpha)$.
Notice that these matrices have size $n+1$, whereas those of section
\ref{hypersurfaces} have size $n+2$.

We now assume that $\bar\cV=\h^n\times\R$ and that $\cV$ is a
Riemannian manifold of dimension $n$ endowed with $\rmS$, $T$, $\nu$
satisfying the compatibility equations for $\h^n\times\R$. We consider
a local orthonormal frame $(e_1,\dots,e_n)$ on $U\subset\cV$, the
associated one-forms $\omega^\alpha$, $\omega^\alpha_\beta$ and the matrix
of one-forms $\Omega\in\cM_{n+1}(\R)$. 

We use the fact that there exists an orthonormal frame on $\h^n$ whose
Christoffel symbols are constant; more precisely, we can choose the
frame $(E_\alpha)$ on $\h^n\times\R$ such that
$\bar\Gamma_{ij}^i=-\bar\Gamma_{ii}^j
=\frac{1}{\sqrt{n}}$ for $i\neq j$,
$i,j\leqslant n$ and all the other Christoffel symbols vanish.

The first step is to prove the following proposition, which is
analogous to proposition \ref{matrixA}.

\begin{prop}
Let $y_0\in\cV$ and $A_0\in\cZ(y_0)$. Then there exist a
neighbourhood $U_1$ of
$y_0$ in $\cV$ and a unique map $A:U_1\to\mathrm{SO}_{n+1}(\R)$ such that
$$A^{-1}\rmd A=\Omega+L(A),$$
$$\forall y\in U_1,\quad A(y)\in\cZ(y),$$
$$A(y_0)=A_0,$$
where $\cZ(y)$ is defined in a way analogous to that of section
\ref{sectiontheorem}.
\end{prop}

To prove this proposition, we introduce the form $\Theta=Z^{-1}\rmd
Z-\Omega-L(Z)$ on 
$\cF=\{(y,Z)\in U\times\mathrm{SO}_{n+1}(\R);Z\in\cZ(y)\}$; this is
well defined since the Christoffel symbols are constant. A long
calculation shows that the distribution $\cD(y,Z)=\ker\Theta_{(y,Z)}$
is involutive. We conclude as in the proof of proposition
\ref{matrixA}.

The second step is to prove the following proposition.

\begin{prop}
Let $x_0\in\h^n\times\R$. There exist a
neighbourhood $U_2$ of $y_0$ contained in $U_1$ and a function
$f:U_2\to\h^n\times\R$ such that
$$\rmd f=(B\circ f)A\omega,$$
$$f(y_0)=x_0,$$
where $\omega$ is the column $(\omega^1,\dots,\omega^n,0)$ and,
for $x\in\h^n\times\R$, $B(x)\in\cM_{n+1}(\R)$ is the matrix
of the coordinates of the frame 
$(E_\alpha(x))$ in the frame
$(\frac{\partial}{\partial x^\alpha})$ (we choose the
upper half-space model for $\h^n$).
\end{prop}

To prove it, we consider the form $B^{-1}\rmd x-A\omega$ on
$U_1\times\bar\cV$, and we show that its kernel again defines an
involutive distribution.

The last step is to check that this map $f$ satisfies the conclusions
of theorem \ref{isometry}.

\section{Applications to minimal surfaces in $\mathbb{M}^2\times\R$}

\subsection{The associate family}

Let $\mathbb{M}^2=\s^2$ or $\mathbb{M}^2=\h^2$. Let $\Sigma$ be a
Riemann surface with a metric $\rmd s^2$ (which we also denote by
$\langle\cdot,\cdot\rangle$), 
$\nabla$ its Riemannian connection, and $\rmJ$ the rotation of angle
$\frac{\pi}{2}$ on $\rmT\Sigma$. Let $\rmS$ be a field of symmetric
operators $\rmS_y:\rmT_y\Sigma\to\rmT_y\Sigma$. Let
$T$ be a vector field on $\Sigma$ 
and $\nu$ a smooth function on $\Sigma$ such that $||T||^2+\nu^2=1$.

\begin{prop} \label{associate}
Assume that $\rmS$ is trace-free and that
$(\rmd s^2,\rmS,T,\nu)$ satisfies the compatibility
equations for $\M^2\times\R$. For $\theta\in\R$ we set
$$\rmS_\theta=e^{\theta\rmJ}\rmS
=(\cos\theta)\rmS+(\sin\theta)\rmJ\rmS,$$
$$T_\theta=e^{\theta\rmJ} T=(\cos\theta)T+(\sin\theta)\rmJ T,$$
i.e., $\rmS_\theta$ and $T_\theta$ are obtained by rotating $\rmS$ and
$T$ by the angle $\theta$.

Then $\rmS_\theta$ is symmetric and trace-free,
$||T_\theta||^2+\nu^2=1$ and 
$(\rmd s^2,\rmS_\theta,T_\theta,\nu)$ satisfies the
compatibility equations for $\M^2\times\R$.
\end{prop}

\begin{proof}
The fact that $\rmS_\theta$ is symmetric and trace-free comes from an
elementary computation. Moreover we have $||T_\theta||=||T||$.
We notice that, since $\dim\Sigma=2$, the Gauss equation
\eqref{gaussmtimesr} is equivalent to 
$$K=\det\rmS+\kappa(1-||T||^2)$$ where $K$ is the Gauss curvature of
$\rmd s^2$. Since $\det(e^{\theta\rmJ})=1$, we have
$\det\rmS_\theta=\det\rmS$, and so the Gauss equation is satisfied for
$(\rmd s^2,\rmS_\theta,T_\theta,\nu)$.

Since $e^{\theta\rmJ}$ commutes with $\nabla_X$ (see \cite{abresch},
section 3.2) 
and  preserves the metric, equations \eqref{conditionT1} and
\eqref{conditionT2} are also satisfied for
$(\rmd s^2,\rmS_\theta,T_\theta,\nu)$.

To prove that the Codazzi equation \eqref{codazzimtimesr} is satisfied
by $(\rmd s^2,\rmS_\theta,T_\theta,\nu)$, we first notice that,
since $\nabla_Xe^{\theta\rmJ}\rmS Y-\nabla_Ye^{\theta\rmJ}\rmS X
-e^{\theta\rmJ}\rmS[X,Y]
=e^{\theta\rmJ}(\nabla_X\rmS Y-\nabla_Y\rmS X-\rmS[X,Y])$, it
suffices to prove that
$$\langle e^{\theta\rmJ}T,Y\rangle X-\langle e^{\theta\rmJ}T,X\rangle Y
=e^{\theta\rmJ}(\langle T,Y\rangle X-\langle T,X\rangle Y).$$
This is obvious at a point where $X=0$. At a point where $X\neq 0$, we can
write $Y=\lambda X+\mu\rmJ X$, and a computation shows that both
expressions are equal to
$\mu\cos\theta\langle T,\rmJ X\rangle X
+\mu\sin\theta\langle T,X\rangle X
-\mu\cos\theta\langle T,X\rangle\rmJ X
+\mu\sin\theta\langle T,\rmJ X\rangle\rmJ X$.
\end{proof}

\begin{thm} \label{associatefamily}
Let $\Sigma$ be a simply connected Riemann surface and
$x:\Sigma\to\M^2\times\R$ a conformal minimal immersion. Let $N$
be the induced normal. Let
$\rmS$ be the symmetric operator on $\Sigma$ induced by the shape
operator of $x(\Sigma)$. Let $T$
be the vector field on $\Sigma$ such that $\rmd x(T)$ is the
projection of $\frac{\partial}{\partial t}$ onto $\rmT(x(\Sigma))$. Let
$\nu=\left\langle N,\frac{\partial}{\partial t}\right\rangle$.

Let $z_0\in\Sigma$.
Then there exists a unique family $(x_\theta)_{\theta\in\R}$ of
conformal minimal immersions 
$x_\theta:\Sigma\to\M^2\times\R$ such that:
\begin{enumerate}
\item $x_\theta(z_0)=x(z_0)$ and $(\rmd x_\theta)_{z_0}=(\rmd
x)_{z_0}$, 
\item the metrics induced on $\Sigma$ by $x$ and $x_\theta$ are the same,
\item the symmetric operator on $\Sigma$ induced by the shape operator
of $x_\theta(\Sigma)$ is $e^{\theta\rmJ}\rmS$,
\item $\frac{\partial}{\partial t}=\rmd x_\theta(e^{\theta\rmJ}T)+\nu
N_\theta$ where $N_\theta$ is the unit normal to $x_\theta$.
\end{enumerate}

Moreover we have $x_0=x$ and the family $(x_\theta)$ is continuous
with respect to $\theta$.

The family of immersions $(x_\theta)_{\theta\in\R}$ is called the
associate family of the immersion $x$. The immersion
$x_{\frac{\pi}{2}}$ is called the conjugate immersion of the immersion
$x$. The immersion
$x_{\pi}$ is called the opposite immersion of the immersion $x$. 
\end{thm}

\begin{proof}
Let $\rmd s^2$ be the metric on $\Sigma$ induced by $x$. Then $(\rmd
s^2,\rmS,T,\nu)$ satisfies the compatibility equations for
$\M^2\times\R$. Thus, by proposition \ref{associate}, 
$(\rmd s^2,e^{\theta\rmJ}\rmS,e^{\theta\rmJ}T,\nu)$ also does. Thus by
theorem \ref{isometry} there exists a unique immersion $x_\theta$
satisfying the properties of the theorem. The fact that $x_0=x$ is
clear.

Finally, $(\rmd s^2,e^{\theta\rmJ}\rmS,e^{\theta\rmJ}T,\nu)$ defines a
matrix of one-forms $\Omega_\theta$ and a matrix of functions
$A_\theta$ satisfying $A_\theta^{-1}\rmd A_\theta=\Omega_\theta$ (by
proposition \ref{matrixA}). By continuity of $\Omega_\theta$ with
respect to $\theta$ we obtain the continuity of $A_\theta$ with
respect to $\theta$, and then the continuity of $x_\theta$ with
respect to $\theta$. 
\end{proof}

\begin{rem} Let $\tau:\Sigma'\to\Sigma$ be a conformal
diffeomorphism. If $\tau$ preserves the orientation, then
$(x\circ\tau)_\theta=x_\theta\circ\tau$; if $\tau$ reverses the
orientation, then $(x\circ\tau)_\theta=x_{-\theta}\circ\tau$.
\end{rem}

In the sequel, we will speak of associate and conjugate immersions
even if condition 1 is not satisfied, i.e., we will consider these
notions up to isometries of $\M^2\times\R$ preserving the orientations
of both $\M^2$ and $\R$.

\begin{rem} The opposite immersion is $x_\pi=\sigma\circ x$ where
$\sigma$ is an isometry of $\M^2\times\R$ preserving the orientation
of $\M^2$ and reversing the orientation of $\R$ (see proposition
\ref{changeofsigns}, case 2).
\end{rem}

\begin{rem}
This associate family for minimal immersions in $\M^2\times\R$ is
analogous to the associate family for minimal immersions in $\R^3$.
Conformal minimal immersions in $\R^3$ are given by the Weierstrass
representation:
$$x(z)=x(z_0)+\re\int_{z_0}^z(1-g^2,i(1+g^2),2g)\omega$$ where $g$ is
a meromorphic function on $\Sigma$ (the Gauss map) and $\omega$ a
holomorphic one-form. Then the associate immersions are
$$x_\theta(z)=x(z_0)
+\re\int_{z_0}^z(1-g^2,i(1+g^2),2g)e^{-i\theta}\omega.$$
\end{rem}

Let $x=(\varphi,h):\Sigma\to\M^2\times\R$ be a conformal minimal
immersion. Then 
$h$ is a real harmonic function and $\varphi$ is a harmonic map to
$\M^2$. We set
$$\frac{\partial}{\partial z}
=\frac 12\left(\frac\partial{\partial u}
-i\frac\partial{\partial v}\right),\quad
\frac{\partial}{\partial\bar z}
=\frac 12\left(\frac\partial{\partial u}
+i\frac\partial{\partial v}\right).$$
The Hopf differential of $\varphi$ is the following $2$-form (see
\cite{rosenbergillinois}): 
$$\rmQ\varphi=4\left\langle\frac{\partial\varphi}{\partial z},
\frac{\partial\varphi}{\partial z}\right\rangle\rmd z^2
=\left(
\left|\left|\frac{\partial\varphi}{\partial u}\right|\right|^2
-\left|\left|\frac{\partial\varphi}{\partial v}\right|\right|^2
-2i\left\langle\frac{\partial\varphi}{\partial u},
\frac{\partial\varphi}{\partial v}\right\rangle\right)\rmd z^2.$$
It is a holomorphic $2$-form on $\Sigma$, and since $x$ is conformal
we have
$$\rmQ\varphi=-4\left(\frac{\partial h}{\partial z}\right)^2\rmd z^2
=-(\rmd(h+ih^*))^2
=-4\left\langle T,\frac{\partial x}{\partial z}\right\rangle\rmd z^2$$
where $h^*$ is the harmonic conjugate function of $h$ (i.e., 
$\frac{\partial h^*}{\partial u}=-\frac{\partial h}{\partial v}$ and
$\frac{\partial h^*}{\partial v}=\frac{\partial h}{\partial u}$). The
reader can refer to \cite{schoenyau} for harmonic maps.

\begin{prop}
Let $x=(\varphi,h):\Sigma\to\M^2\times\R$ be a conformal minimal immersion,
and $(x_\theta)=((\varphi_\theta,h_\theta))$ its associate family of
conformal minimal immersions. Let $h^*$ be the harmonic conjugate of
$h$. Then we have
$$h_\theta=(\cos\theta)h+(\sin\theta)h^*,\quad
\rmQ\varphi_\theta=e^{-2i\theta}\rmQ\varphi.$$
\end{prop}

\begin{proof}
We have 
\begin{eqnarray*}
\frac{\partial h_\theta}{\partial u}
=\left\langle\frac{\partial x_\theta}{\partial u},
\frac{\partial}{\partial t}\right\rangle
=\left\langle\frac{\partial}{\partial u},T_\theta\right\rangle 
& = & \cos\theta\left\langle\frac{\partial}{\partial u},T\right\rangle
+\sin\theta\left\langle\frac{\partial}{\partial u},\rmJ T\right\rangle
\\
& = & \cos\theta\left\langle\frac{\partial}{\partial u},T\right\rangle
-\sin\theta\left\langle\frac{\partial}{\partial v},T\right\rangle \\
& = & \cos\theta\frac{\partial h}{\partial u}
-\sin\theta\frac{\partial h}{\partial v}.
\end{eqnarray*}
In the same way we have $\frac{\partial h_\theta}{\partial v}
=\cos\theta\left\langle\frac{\partial}{\partial v},T\right\rangle
+\sin\theta\left\langle\frac{\partial}{\partial v},\rmJ T\right\rangle
=\cos\theta\frac{\partial h}{\partial v}
+\sin\theta\frac{\partial h}{\partial u}$. This proves that
$h_\theta=(\cos\theta)h+(\sin\theta)h^*$. The expression of
$\rmQ\varphi_\theta$ follows immediately. 
\end{proof}

\begin{rem}
Recently, Hauswirth, S\'a Earp and Toubiana (\cite{hset}) defined the
following notion of associated immersions in $\h^2\times\R$: two
isometric conformal minimal immersions in $\h^2\times\R$ are said to
be associated if their Hopf differential differ by the multiplication
by some constant $e^{i\theta}$. Morover, they proved that two
isometric conformal minimal immersions in $\h^2\times\R$ having the
same Hopf differential are equal up to an isometry of $\h^2\times\R$.
Thus the notions of associated immersions in the sense of this paper
and in the sense of \cite{hset} are equivalent.

In \cite{screw}, S\'a Earp and Toubiana ask the following question: if
two conformal minimal immersions $x,\tilde x:\Sigma\to\M^2\times\R$
are isometric, then are they associated ? (This result holds for $\R^3$.)
\end{rem}

\begin{rem}
Abresch and Rosenberg (\cite{abresch}) defined a holomorphic Hopf
differential for constant mean curvature surfaces in
$\M^2\times\R$. For minimal surfaces 
in $\M^2\times\R$, this Hopf differential is
\begin{eqnarray*}
Q(X,Y) & = &
-\frac{\kappa}{2}(\langle T,X\rangle\langle T,Y\rangle
-\langle T,\rmJ X\rangle\langle T,\rmJ Y\rangle) \\
& & +i\frac{\kappa}{2}(\langle T,\rmJ X\rangle\langle T,Y\rangle
+\langle T,X\rangle\langle T,\rmJ Y\rangle).
\end{eqnarray*}
A computation shows that $$Q=\frac{\kappa}{2}\rmQ\varphi.$$
\end{rem}

\begin{prop} \label{zeroT}
Let $x:\Sigma\to\M^2\times\R$ be a conformal minimal immersion. If $x$
does not define a horizontal $\M^2\times\{t\}$, then the zeros of $T$
are isolated. 
\end{prop}

\begin{proof}
The height function $h=\langle x,\frac\partial{\partial t}\rangle$
satisfies 
$\rmd h(X)=\langle T,X\rangle$; thus the zeroes of $T$ are the zeroes
of $\rmd h$. Since $h$ is harmonic, either the zeroes of $\rmd h$ are
isolated or $h$ is constant. The latter case is excluded by
hypothesis.
\end{proof}

\begin{rem}
Umbilic points (i.e., zeroes of the shape operator) may be
non-isolated: for example helicoids and unduloids in $\s^2\times\R$
have curves of umbilic points (see section \ref{exampless2}).
\end{rem}

We now give some geometric properties of conjugate surfaces.

The transformation $\rmS\mapsto\rmJ\rmS$ implies that curvature lines
and asymptotic lines are exchanged by conjugation (as in $\R^3$). (More
generally the normal curvature and the normal torsion of a curve are
swapped up to a sign.) The reader can refer to \cite{karcher} for
geometric properties of conjugate surfaces in $\R^3$.

Moreover, the transformation $T\mapsto\rmJ T$ implies the following
transformation: a horizontal curve $\gamma$ along which the surface is
vertical (i.e., $\nu=0$ along $\gamma$ and $\gamma'$ is orthogonal to
$T$) is mapped to vertical curve (i.e., $\nu=0$ along $\gamma$ and
$\gamma'$ is
proportional to $T$), and vice versa. We also notice that a minimal
surface cannot be horizontal along a horizontal curve unless the
minimal surface is a horizontal $\M^2\times\{t\}$ (indeed, this would
imply that $T=0$ along this curve).

Hence conjugation swaps two pairs of Schwarz reflections:
\begin{enumerate}
\item the symmetry with respect to a vertical plane containing a
curvature line becomes the rotation with respect to a horizontal
geodesic of $\M^2$,
and vice versa,
\item the symmetry with respect to a horizontal plane containing a
curvature line becomes the rotation with respect to a vertical
straight line, and vice versa.
\end{enumerate}
The first case is illustrated by a generatrix curve of an unduloid or a
catenoid and a horizontal line of a helicoid; the second case is
illustrated by the waist circle of an unduloid or a catenoid and the
axis of a helicoid. These examples are detailed in sections
\ref{exampless2} and \ref{examplesh2}.

\subsection{Helicoids and unduloids in $\s^2\times\R$}
\label{exampless2}

Apart from the horizontal spheres
$\s^2\times\{t\}$ and the vertical cylinders $\s^1\times\R$ ($\s^1$
being a great circle in $\s^2$), the
most simple examples of minimal surfaces in $\s^2\times\R$ are
helicoids and unduloids. Theses surfaces are described in 
\cite{pedrosa} and \cite{rosenbergillinois}. They are properly
embedded and foliated by circles. Unduloids are rotational and
vertically periodic; helicoids are invariant by a screw motion.

\subsubsection*{Helicoids.}
For $\beta\neq 0$, the helicoid $\cH_\beta$ is given by the following
conformal immersion:
$$x(u,v)=\left(\begin{array}{c}
\sin\varphi(u)\cos\beta v \\
\sin\varphi(u)\sin\beta v \\
\cos\varphi(u) \\
v
\end{array}\right),$$
where the function $\varphi$ satisfies
\begin{equation} \label{equationphi}
\varphi'(u)^2=1+\beta^2\sin^2\varphi(u),\quad
\varphi''(u)=\beta^2\sin\varphi(u)\cos\varphi(u).
\end{equation}
We can assume that $\varphi(0)=0$ and $\varphi'(u)>0$. When $\beta>0$
we say that $\cH_\beta$ is a right helicoid; when $\beta<0$
we say that $\cH_\beta$ is a left helicoid.

The normal to $\s^2\times\R$ in $\R^4$ is
$$\bar N(u,v)=\left(\begin{array}{c}
\sin\varphi(u)\cos\beta v \\
\sin\varphi(u)\sin\beta v \\
\cos\varphi(u) \\
0
\end{array}\right).$$
The normal to $\cH_\beta$ in $\s^2\times\R$ is
$$N(u,v)=\frac{1}{\varphi'(u)}\left(\begin{array}{c}
\sin\beta v \\
-\cos\beta v \\
0 \\
\beta\sin\varphi(u)
\end{array}\right).$$

We compute:
$$\left\langle\frac{\partial^2x}{\partial u^2},N\right\rangle
=\left\langle\frac{\partial^2x}{\partial v^2},N\right\rangle=0,\quad
\left\langle\frac{\partial^2x}{\partial u\partial v},N\right\rangle
=-\beta\cos\varphi(u).$$ Using the fact that
$\langle\rmS X,Y\rangle=\langle\rmd Y(X),N\rangle$,
we compute that the matrix of $\rmS$ in the frame
$(\frac{\partial}{\partial u},\frac{\partial}{\partial v})$ is the
following: 
$$-\frac{\beta\cos\varphi(u)}{\varphi'(u)^2}
\left(\begin{array}{cc}
0 & 1 \\
1 & 0
\end{array}\right).$$
In particular the points where $\cos\varphi(u)=0$ are umbilic points.
We also have $$T=\frac{1}{\varphi'(u)^2}\frac{\partial}{\partial v},
\quad\nu=\frac{\beta\sin\varphi(u)}{\varphi'(u)}.$$ 

\begin{rem}
When $\beta=0$, the formula defines a vertical cylinder
$\s^1\times\R$. When $\beta\to\infty$, the surface converges to the
foliation by horizontal spheres $\s^2\times\{t\}$. 
\end{rem}

\subsubsection*{Unduloids.}
For $\alpha>1$ or $\alpha<-1$, the unduloid $\cU_\alpha$ is given by
the following conformal immersion:
$$x(u,v)=\left(\begin{array}{c}
\sin\psi(u)\cos\alpha v \\
\sin\psi(u)\sin\alpha v \\
\cos\psi(u) \\
u
\end{array}\right),$$
where the function $\psi$ satisfies
\begin{equation} \label{equationpsi}
1+\psi'(u)^2=\alpha^2\sin^2\psi(u),\quad
\psi''(u)=\alpha^2\sin\psi(u)\cos\psi(u).
\end{equation}
We can assume that $\psi'(0)=0$, $\psi(u)\in(0,\pi)$ and
$\cos\psi(0)>0$.

The normal to $\cU_\alpha$ in $\s^2\times\R$ is
$$N(u,v)=\frac{1}{\alpha\sin\psi(u)}\left(\begin{array}{c}
-\cos\psi(u)\cos\alpha v \\
-\cos\psi(u)\sin\alpha v \\
\sin\psi(u) \\
\psi'(u)
\end{array}\right).$$
We compute that the matrix of $\rmS$ in the frame
$(\frac{\partial}{\partial u},\frac{\partial}{\partial v})$ is the
following: 
$$-\frac{\alpha\cos\psi(u)}{1+\psi'(u)^2}
\left(\begin{array}{cc}
1 & 0 \\
0 & -1
\end{array}\right).$$
In particular the points where $\cos\psi(u)=0$ are umbilic points.
We also have $$T=\frac{1}{1+\psi'(u)^2}\frac{\partial}{\partial u},
\quad\nu=\frac{\psi'(u)}{\alpha\sin\psi(u)}.$$

\begin{rem}
When $\alpha=\pm 1$, the formula defines a vertical cylinder
$\s^1\times\R$. When $\alpha\to\infty$, the surface converges to the
foliation by horizontal spheres $\s^2\times\{t\}$. 
\end{rem}

\begin{prop}
The conjugate surface of the unduloid $\cU_\alpha$ is the helicoid
$\cH_\beta$ with $\alpha^2=1+\beta^2$ and $\alpha$, $\beta$ having the
same sign.
\end{prop}

\begin{proof}
We set $y_1(u)=\alpha\cos\psi(u)$ and
$y_2(u)=\beta\cos\varphi(u)$. A computation shows that both $y_1$ and
$y_2$ are solutions of the following equation:
$$(y')^2=(y^2-\alpha^2)(y^2-\beta^2),$$
and hence of the following equation:
$$y''=y(2y^2-\alpha^2-\beta^2).$$
We have $\psi'(0)=0$ and so by \eqref{equationpsi} we have
$y_1(0)^2=\beta^2$ and thus $y_1'(0)=0$, and $\varphi(0)=0$ so
$y_2(0)=\beta$ and thus $y_2'(0)=0$. Moreover, $\cos\psi(0)>0$, so
$y_1(0)$ has the sign of $\alpha$; since $\alpha$ and $\beta$ have the
same sign, we have $y_1(0)=\beta$. By the Cauchy-Lipschitz theorem 
we conclude that $y_1=y_2$. From this we deduce using
\eqref{equationpsi} and \eqref{equationphi} that
$\varphi'(u)^2=1+\psi'(u)^2$, and thus 
$\cU_\alpha$ and $\cH_\beta$ are locally isometric, and
$\rmS_{\cH_\beta}=\rmJ\rmS_{\cU_\alpha}$ and
$T_{\cH_\beta}=\rmJ T_{\cU_\alpha}$.
Finally we have $\nu_{\cU_\alpha}=-\frac{y_1'}{\alpha^2-y_1^2}$
and $\nu_{\cH_\beta}=-\frac{y_2'}{\alpha^2-y_2^2}$,  
so we get $\nu_{\cH_\beta}=\nu_{\cU_\alpha}$.
\end{proof}

\begin{rem}
The vertical cylinder $\s^1\times\R$ is globally invariant by
conjugation, but the vertical lines and the horizontal circles are
exchanged. For example, a rectangle of height $t$ and whose basis is
an arc of angle $\theta$ becomes a rectangle of height $\theta$ and
whose basis is an arc of angle $t$. 

The horizontal sphere $\s^2\times\{0\}$ is pointwise invariant by
conjugation (since it satisfies $\rmS=0$ and $T=0$).
\end{rem}

\begin{rem}
The horizontal projections of helicoids and unduloids are the Gauss
maps of constant mean curvature Delaunay surfaces in $\R^3$:
helicoids in $\s^2\times\R$ come from nodoids in $\R^3$ and unduloids
in $\s^2\times\R$ come from unduloids in $\R^3$. This correspondance
is described in \cite{rosenbergimpa}.
\end{rem}

\subsection{Helicoids and generalized catenoids in $\h^2\times\R$}
\label{examplesh2}

Apart from the horizontal planes
$\h^2\times\{t\}$ and the vertical planes $\h^1\times\R$ ($\h^1$
being a geodesic of $\h^2$), the
most simple examples of minimal surfaces in $\h^2\times\R$ are
helicoids and catenoids. These surfaces are described in
\cite{pedrosa} and \cite{nelli}. They are properly embedded. Catenoids
are rotational; helicoids are invariant by a screw motion and foliated
by geodesics of $\h^2$.

More generally, Hauswirth classified minimal surfaces
in $\h^2\times\R$ foliated by horizontal curves of constant curvature
in $\h^2$ (\cite{hauswirth}). These surfaces form a two-parameter
family. This family includes, among others, catenoids, helicoids and
Riemann-type examples. All the surfaces described in this section
belong to the Hauswirth family.

\subsubsection*{Helicoids.}
For $\beta\neq 0$, the helicoid $\cH_\beta$ is given by the following
conformal immersion:
$$x(u,v)=\left(\begin{array}{c}
\cosh\varphi(u) \\
\sinh\varphi(u)\cos\beta v \\
\sinh\varphi(u)\sin\beta v \\
v
\end{array}\right),$$
where the function $\varphi$ satisfies
\begin{equation} \label{equationphihn}
\varphi'(u)^2=1+\beta^2\sinh^2\varphi(u),\quad
\varphi''(u)=\beta^2\sinh\varphi(u)\cosh\varphi(u).
\end{equation}
We can assume that $\varphi(0)=0$ and $\varphi'(u)>0$. The function
$\varphi$ is defined on a bounded interval. When $\beta>0$ we
say that $\cH_\beta$ is a right 
helicoid; when $\beta<0$ we say that $\cH_\beta$ is a left
helicoid.

The normal to $\cH_\beta$ in $\h^2\times\R$ is
$$N(u,v)=\frac{1}{\varphi'(u)}\left(\begin{array}{c}
0 \\
\sin\beta v \\
-\cos\beta v \\
\beta\sinh\varphi(u)
\end{array}\right).$$
now $\beta>0$.
we compute that the matrix of $\rmS$ in the frame
$(\frac{\partial}{\partial u},\frac{\partial}{\partial v})$ is the
following: 
$$-\frac{\beta\cosh\varphi(u)}{\varphi'(u)^2}
\left(\begin{array}{cc}
0 & 1 \\
1 & 0
\end{array}\right).$$
We also have $$T=\frac{1}{\varphi'(u)^2}\frac{\partial}{\partial v},
\quad\nu=\frac{\beta\sinh\varphi(u)}{\varphi'(u)}.$$ 

\begin{rem}
When $\beta=0$, the fomula defines a vertical plane
$\h^1\times\R$. When $\beta\to\infty$, the surface converges to the 
foliation by horizontal planes $\h^2\times\{t\}$. 
\end{rem}

\subsubsection*{Catenoids.}
For $\alpha\neq 0$, the catenoid $\cC_\alpha$ is given by the
following conformal immersion:
$$x(u,v)=\left(\begin{array}{c}
\cosh\psi(u) \\
\sinh\psi(u)\cos\alpha v \\
\sinh\psi(u)\sin\alpha v \\
u
\end{array}\right),$$
where the function $\psi$ satisfies
\begin{equation} \label{equationpsihn}
1+\psi'(u)^2=\alpha^2\sinh^2\psi(u),\quad
\psi''(u)=\alpha^2\sinh\psi(u)\cosh\psi(u).
\end{equation}
We can assume that $\psi'(0)=0$ and $\psi(u)>0$. The function
$\psi$ is defined on the interval $(-u_0,u_0)$
with
$$u_0=\int_{\psi(0)}^\infty\frac{\rmd\psi}
{\sqrt{\alpha^2\sinh^2\psi-1}}
=\int_1^\infty\frac{\rmd x}{\sqrt{(x^2+\alpha^2)(x^2-1)}}.$$
Thus we have $$u_0<\int_1^\infty\frac{\rmd x}{x\sqrt{x^2-1}}
=\frac{\pi}{2}.$$
This proves that the height of the catenoid $\cC_\alpha$ is smaller
than $\pi$; moreover the height tends to $0$ when $\alpha\to\infty$
and to $\pi$ when $\alpha\to 0$ (theorem 1 in \cite{nelli} holds for
$t\in(0,\frac{\pi}{2})$).
The function $\psi$ is decreasing on $(-u_0,0)$ and increasing on
$(0,u_0)$. The waist circle is given by $u=0$. 

The normal to $\cC_\alpha$ in $\h^2\times\R$ is
$$N(u,v)=\frac{1}{\alpha\sinh\psi(u)}\left(\begin{array}{c}
-\sinh\psi(u) \\
-\cosh\psi(u)\cos\alpha v \\
-\cosh\psi(u)\sin\alpha v \\
\psi'(u)
\end{array}\right).$$
We compute that the matrix of $\rmS$ in the frame
$(\frac{\partial}{\partial u},\frac{\partial}{\partial v})$ is the
following: 
$$-\frac{\alpha\cosh\psi(u)}{1+\psi'(u)^2}
\left(\begin{array}{cc}
1 & 0 \\
0 & -1
\end{array}\right).$$
We also have
$$T=\frac{1}{1+\psi'(u)^2}\frac{\partial}{\partial u},\quad
\nu=\frac{\psi'(u)}{\alpha\sinh\psi(u)}.$$

\subsubsection*{A minimal surface foliated by horocycles.}
We search a minimal surface such that each horizontal curve is a
horocycle in $\h^2$ and such that all the horocycles have the same
asymptotic 
point. Such a surface can be parametrized in the following way:
$$x(u,v)=\left(\begin{array}{c}
\frac{\lambda(u)}{2}+\frac{1+f(u,v)^2}{2\lambda(u)} \\
f(u,v) \\
-\frac{\lambda(u)}{2}+\frac{1+f(u,v)^2}{2\lambda(u)} \\
u
\end{array}\right)$$
with $\lambda>0$ and $\frac{\partial f}{\partial v}>0$.
This immersion is conformal if and only if
$$\frac{\partial f}{\partial u}=\frac{f\lambda'}{\lambda},
\quad\left(\frac{\partial f}{\partial v}\right)^2
=1+\left(\frac{\lambda'}{\lambda}\right)^2.$$
We deduce from the second relation that
$\frac{\partial^2f}{\partial v^2}=0$, and so
$$f(u,v)=\alpha(u)v+\beta(u).$$
Reporting in the first relation we get
$$\frac{\alpha'}{\alpha}=\frac{\beta'}{\beta}
=\frac{\lambda'}{\lambda}.$$
The immersion is minimal if and only if $\Delta x$ is proportional to
the normal $\bar N$ to $\h^2\times\R$; a computation shows that this
happens if and only if
$(\lambda')^2+\alpha^2\lambda^2=\lambda\lambda''$, i.e., if and only
if $2(\lambda')^2+\lambda^2=\lambda\lambda''$, or, equivalently,
$$\left(\frac{1}{\lambda}\right)''=-\frac{1}{\lambda}.$$
Up to a reparametrization and an isometry of $\h^2$ we can choose
$\lambda(u)=\alpha(u)=\frac{1}{\cos u}$ for
$u\in(-\frac{\pi}{2},\frac{\pi}{2})$ and $\beta(u)=0$. Thus we get the
following proposition.

\begin{prop}
The map
$$x(u,v)=\left(\begin{array}{c}
\frac{v^2+1}{2\cos u}+\frac{\cos u}{2} \\
\frac{v}{\cos u} \\
\frac{v^2-1}{2\cos u}+\frac{\cos u}{2} \\
u
\end{array}\right)$$
defined for $(u,v)\in(-\frac{\pi}{2},\frac{\pi}{2})\times\R$ 
is a conformal
minimal embedding such that the curves $u=u_0$ are horocycles in
$\h^2$ having the same asymptotic point. We will denote this surface
by $\cC_0$.

Morover, the surface $\cC_0$ is the unique one (up to isometries of
$\h^2\times\R$) having this property.
\end{prop}

\begin{figure}[htbp] 
\begin{center}
\includegraphics{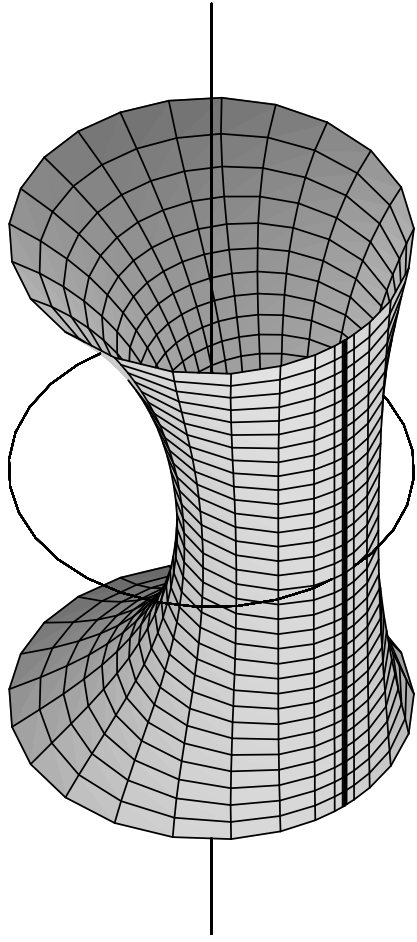}
\caption{A minimal surface in $\h^2\times\R$ foliated by horocycles.}
\label{horocycles}
\end{center}
\end{figure}

In the upper half-plane model for $\h^2$, the curve at height $u$ of
$\cC_0$ is the horizontal Euclidean line $x_2=\cos u$. Figure
\ref{horocycles} is a picture of $\cC_0$ (in this picture the model
for $\h^2$ is the Poincar\'e unit disk model). The surface $\cC_0$ has
height $\pi$. It is symmetric with respect to the horizontal plane
$\h^2\times\{0\}$ and it is invariant by a one-parameter family of
horizontal parabolic isometries.

The normal to $\cC_0$ in $\h^2\times\R$ is
$$N(u,v)=\left(\begin{array}{c}
-\frac{v^2+1}{2}+\frac{\cos^2u}{2} \\
-v \\
\frac{1-v^2}{2}+\frac{\cos^2u}{2} \\
\sin u
\end{array}\right).$$
We compute that the matrix of $\rmS$ in the frame
$(\frac{\partial}{\partial u},\frac{\partial}{\partial v})$ is the
following: 
$$-\cos u
\left(\begin{array}{cc}
1 & 0 \\
0 & -1
\end{array}\right).$$
We also have $$T=\cos^2u\frac{\partial}{\partial u},
\quad\nu=\sin u.$$ 

\subsubsection*{Minimal surfaces foliated by equidistants.}
For $\gamma\in(0,1)$ or $\gamma\in(-1,0)$, we consider the following
immersion: 
$$x(u,v)=\left(\begin{array}{c}
\cosh\chi(u)\cosh\gamma v \\
\sinh\chi(u) \\
\cosh\chi(u)\sinh\gamma v \\
u
\end{array}\right).$$
with
\begin{equation} \label{equationchi}
1+\chi'(u)^2=\gamma^2\cosh^2\chi(u),\quad
\chi''(u)=\gamma^2\cosh\chi(u)\sinh\chi(u).
\end{equation}
It is a conformal minimal immersion. 

We choose $\chi$ such that $\chi'(0)=0$ and $\chi(u)>0$. 
The function
$\chi$ is defined on the interval $(-u_0,u_0)$
with
$$u_0=\int_{\chi(0)}^\infty\frac{\rmd\chi}
{\sqrt{\gamma^2\cosh^2\chi-1}}
=\int_1^\infty\frac{\rmd x}{\sqrt{(x^2-\gamma^2)(x^2-1)}}.$$
Thus we have $$u_0>\int_1^\infty\frac{\rmd x}{x\sqrt{x^2-1}}
=\frac{\pi}{2}.$$
We have defined a minimal surface $\cG_\gamma$, which we call a
generalized catenoid. Its height is greater
than $\pi$, tends to $\pi$ when $\gamma\to 0$
and to $+\infty$ when $\gamma\to 1$. 
The function $\chi$ is decreasing on $(-u_0,0)$ and increasing on
$(0,u_0)$. The surface is symmetric with respect to the horizontal
plane $\h^2\times\{0\}$ and it is invariant by a one-parameter family
of horizontal hyperbolic isometries. The horizontal curves are
equidistants to a geodesic in $\h^2$.

The normal to $\cG_\gamma$ in $\h^2\times\R$ is
$$N(u,v)=-\frac{1}{\gamma\cosh\chi(u)}\left(\begin{array}{c}
\sinh\chi(u)\cosh\gamma v \\
\cosh\chi(u) \\
\sinh\chi(u)\sinh\gamma v \\
-\chi'(u)
\end{array}\right).$$
We compute that the matrix of $\rmS$ in the frame
$(\frac{\partial}{\partial u},\frac{\partial}{\partial v})$ is the
following: 
$$-\frac{\gamma\sinh\chi(u)}{1+\chi'(u)^2}
\left(\begin{array}{cc}
1 & 0 \\
0 & -1
\end{array}\right).$$
We also have $$T=\frac{1}{1+\chi'(u)^2}\frac{\partial}{\partial u},
\quad\nu=\frac{\chi'(u)}{\gamma\cosh\chi(u)}.$$ 

\begin{rem}
When $\gamma=\pm 1$, the formula defines a vertical plane
$\h^1\times\R$.
\end{rem}

\begin{prop}
The conjugate surface of the catenoid $\cC_\alpha$ is the helicoid
$\cH_\beta$ with $\beta^2=1+\alpha^2$ and $\alpha$, $\beta$ having the
same sign.
\end{prop}

\begin{proof}
We set $y_1(u)=\alpha\cosh\psi(u)$ and
$y_2(u)=\beta\cosh\varphi(u)$. A computation shows that both $y_1$ and
$y_2$ are solutions of the following equation:
$$(y')^2=(y^2-\alpha^2)(y^2-\beta^2),$$
and hence of the following equation:
$$y''=y(2y^2-\alpha^2-\beta^2).$$
We have $\psi'(0)=0$ and so by \eqref{equationpsihn} we have
$y_1(0)^2=\beta^2$ and thus $y_1'(0)=0$, and $\varphi(0)=0$ so
$y_2(0)=\beta$ and thus $y_2'(0)=0$. Moreover, $y_1(0)$ has the sign
of $\alpha$, i.e., the sign of $\beta$, so we get $y_1(0)=\beta$. By
the Cauchy-Lipschitz theorem 
we conclude that $y_1=y_2$ (and in particular they have the same
domain of definition). From this we deduce using \eqref{equationpsihn}
and \eqref{equationphihn} that $\varphi'(u)^2=1+\psi'(u)^2$, and thus
$\cC_\alpha$ and $\cH_\beta$ are locally isometric,
$\rmS_{\cH_\beta}=\rmJ\rmS_{\cC_\alpha}$ and $T_{\cH_\beta}=\rmJ
T_{\cC_\alpha}$. Finally we have
$\nu_{\cC_\alpha}=\frac{y_1'}{y_1^2-\alpha^2}$ and
$\nu_{\cH_\beta}=\frac{y_2'}{y_2^2-\alpha^2}$, so we get
$\nu_{\cH_\beta}=\nu_{\cC_\alpha}$.
\end{proof}

\begin{prop}
The conjugate surface of the surface $\cC_0$ is the helicoid $\cH_1$.
\end{prop}

\begin{proof}
In the case where $\beta=1$, the function $\varphi$ satisfies
$\varphi'=\cosh\varphi$, and thus we have
$\varphi(u)=\ln(\tan(\frac{u}{2}+\frac{\pi}{4}))$,
$\varphi'(u)=\frac{1}{\cos u}$ and $\sinh\varphi(u)=\tan u$. Then,
using the above calculations, we
easily check that $\cC_0$ and $\cH_1$ are locally isometric, and that
$\rmS_{\cH_1}=\rmJ\rmS_{\cC_0}$, $T_{\cH_1}=\rmJ T_{\cC_0}$,
$\nu_{\cH_1}=\nu_{\cC_0}$.
\end{proof} 

\begin{rem}
The conjugate surface of the surface $\cC_0$ with the opposite
orientation is the helicoid $\cH_{-1}$.
\end{rem}

\begin{prop}
The conjugate surface of the generalized catenoid $\cG_\gamma$ is the
helicoid 
$\cH_\beta$ with $\beta^2+\gamma^2=1$ and $\beta$, $\gamma$ having the
same sign.
\end{prop}

\begin{proof}
We set $y_1(u)=\gamma\sinh\chi(u)$ and
$y_2(u)=\beta\cosh\varphi(u)$. A computation shows that both $y_1$ and
$y_2$ are solutions of the following equation:
$$(y')^2=(y^2+\gamma^2)(y^2-\beta^2),$$
and hence of the following equation:
$$y''=y(2y^2+\gamma^2-\beta^2).$$
We have $\chi'(0)=0$ and  so by \eqref{equationchi} we have
$y_1(0)^2=\beta^2$ and thus $y_1'(0)=0$, and $\varphi(0)=0$ so
$y_2(0)=\beta$ and thus $y_2'(0)=0$. Moreover, $y_1(0)$ has the sign
of $\gamma$, i.e., the sign of $\beta$, so we get $y_1(0)=\beta$. By
the Cauchy-Lipschitz theorem 
we conclude that $y_1=y_2$ (and in particular they have the same
domain of definition). From this we deduce using \eqref{equationchi}
and \eqref{equationphihn} that $\varphi'(u)^2=1+\chi'(u)^2$, and thus
$\cG_\gamma$ and $\cH_\beta$ are locally isometric,
$\rmS_{\cH_\beta}=\rmJ\rmS_{\cG_\gamma}$ and
$T_{\cH_\beta}=\rmJ T_{\cG_\gamma}$.
Finally we have
$\nu_{\cG_\gamma}=\frac{y_1'}{y_1^2+\gamma^2}$ and
$\nu_{\cH_\beta}=\frac{y_2'}{y_2^2+\gamma^2}$, so we get
$\nu_{\cH_\beta}=\nu_{\cG_\gamma}$.
\end{proof}

\begin{rem}
This study shows that there are three types of helicoid conjugates
according to the parameter of the screw-motion associated to the
helicoid: the first type ones are the catenoids, which are rotational
surfaces, the second type one is $\cC_0$, which is invariant by a
one-parameter family of
horizontal parabolic isometries and which corresponds to a critical
value of the parameter, the third type ones are the generalized
catenoids, which are invariant a one-parameter family of horizontal
hyperbolic isometries.

This phenomenon is very similar to what happens for the conjugate
cousins in $\h^3$ of the helicoids in $\R^3$. There exists an
isometric correspondance between minimal surfaces in $\R^3$ and
constant mean curvature one surfaces in $\h^3$ 
called the cousin relation (see \cite{bryant} and
\cite{umehara}). Starting from a helicoid in $\R^3$, we 
consider its conjugate surface, which is a catenoid in $\R^3$, and then the
cousin surface in $\h^3$, which is a catenoid cousin. Catenoid cousins
are of three types according to the parameter of the minimal helicoid:
some are rotational surfaces, one is invariant by a one-parameter
family of parabolic isometries (and corresponds to a critical value of
the parameter), some are invariant by a one-parameter
family of hyperbolic isometries. These surfaces are described in
details in \cite{toubiana} and \cite{rosenberg}.
\end{rem}

\begin{rem}
All the above surfaces belong to the Hauswirth family: with the
notations of \cite{hauswirth}, helicoids correspond to $d=0$, $c>0$,
$c\neq 1$; 
catenoids correspond to $c=0$, $d>1$; the surface $\cC_0$ corresponds
to $c=0$, $d=1$; the surfaces $\cG_\gamma$ correspond to $c=0$,
$d\in(0,1)$.
\end{rem}

\begin{rem}
The vertical plane $\h^1\times\R$ is globally invariant by
conjugation, but the vertical lines and the horizontal geodesics of
$\h^2$ are exchanged.
The horizontal plane $\h^2\times\{0\}$ is pointwise invariant by
conjugation (since it satisfies $\rmS=0$ and $T=0$). This is similar
to what happens in $\s^2\times\R$.
\end{rem}

\bibliographystyle{alpha}
\bibliography{codazzi}

\begin{thebibliography}{HSET04}

\bibitem[AR03]{abresch}
U.~Abresch and H.~Rosenberg.
\newblock The {H}opf differential for constant mean curvature surfaces in
  $\mathbb{S}^2\times\mathbb{R}$ and $\mathbb{H}^2\times\mathbb{R}$.
\newblock Preprint, 2003.

\bibitem[Bry87]{bryant}
R.~Bryant.
\newblock Surfaces of mean curvature one in hyperbolic space.
\newblock {\em Astérisque}, 154--155:321--347, 1987.

\bibitem[Car92]{docarmo}
M.~do Carmo.
\newblock {\em Riemannian Geometry}.
\newblock Birkhäuser, 1992.

\bibitem[Hau03]{hauswirth}
L.~Hauswirth.
\newblock Generalized {R}iemann examples in three-dimensional manifolds.
\newblock Preprint, 2003.

\bibitem[HSET04]{hset}
L.~Hauswirth, R.~S\'a~Earp, and \'E. Toubiana.
\newblock Minimal surfaces in $\mathbb{H}^2\times\mathbb{R}$.
\newblock Work in progress, 2004.

\bibitem[Kar01]{karcher}
H.~Karcher.
\newblock Introduction to conjugate {P}lateau constructions.
\newblock Preprint, 2001.

\bibitem[MR03]{meeksrosenberg}
W.~Meeks and H.~Rosenberg.
\newblock The theory of minimal surfaces in ${M}\times\mathbb{R}$.
\newblock Preprint, 2003.

\bibitem[NR02]{nelli}
B.~Nelli and H.~Rosenberg.
\newblock Minimal surfaces in $\mathbb{H}^2\times\mathbb{R}$.
\newblock {\em Bull. Braz. Math. Soc., New Series}, 33(2):263--292, 2002.

\bibitem[PR99]{pedrosa}
R.~Pedrosa and M.~Ritor\'e.
\newblock Isoperimetric domains in the {R}iemannian product of a circle with a
  simply connected space form and applications to free boundary value problems.
\newblock {\em Indiana Univ. Math. J.}, 48:1357--1394, 1999.

\bibitem[Ros02a]{rosenberg}
H.~Rosenberg.
\newblock {B}ryant surfaces.
\newblock In {\em The Global Theory of Minimal Surfaces in Flat Spaces, Martina
  Franca, Italy 1999, Lecture Notes in Mathematics 1775}. Springer, 2002.

\bibitem[Ros02b]{rosenbergillinois}
H.~Rosenberg.
\newblock Minimal surfaces in $\mathbb{M}^2\times\mathbb{R}$.
\newblock {\em Illinois J. Math.}, 46(4):1177--1195, 2002.

\bibitem[Ros03]{rosenbergimpa}
H.~Rosenberg.
\newblock Some recent developments in the theory of minimal surfaces in
  $3$-manifolds.
\newblock In {\em 24$\sp {\rm o}$ Colóquio Brasileiro de Matemática, Instituto
  de Matemática Pura e Aplicada, Rio de Janeiro}. Publicações Matemáticas do
  IMPA, 2003.

\bibitem[SET01]{toubiana}
R.~S\'a~Earp and \'E. Toubiana.
\newblock On the geometry of constant mean curvature one surfaces in hyperbolic
  space.
\newblock {\em Illinois J. Math.}, 45(2):371--401, 2001.

\bibitem[SET04]{screw}
R.~S\'a~Earp and \'E. Toubiana.
\newblock Screw motion surfaces in $\mathbb{H}^2\times\mathbb{R}$ and
  $\mathbb{S}^2\times\mathbb{R}$.
\newblock Preprint, 2004.

\bibitem[SY97]{schoenyau}
R.~Schoen and S.~T. Yau.
\newblock {\em Lectures on Harmonic Maps}.
\newblock Conference Proceedings and Lecture Notes in Geometry and Topology,
  II. International Press, Cambridge, MA, 1997.

\bibitem[Ten71]{tenenblat}
K.~Tenenblat.
\newblock On isometric immersions of {R}iemannian manifolds.
\newblock {\em Bol. Soc. Brasil. Mat.}, 2(2):23--36, 1971.

\bibitem[UY93]{umehara}
M.~Umehara and K.~Yamada.
\newblock Complete surfaces of constant mean curvature $1$ in the hyperbolic
  $3$-space.
\newblock {\em Ann. of Math. (2)}, 137:611--638, 1993.

\end{thebibliography}

\end{document}